# SEMIPARAMETRIC DETECTION OF SIGNIFICANT ACTIVATION FOR BRAIN FMRI

By Chunming Zhang[1] and Tao Yu[2]

*University of Wisconsin-Madison*

Functional magnetic resonance imaging (fMRI) aims to locate activated regions in human brains when specific tasks are performed. The conventional tool for analyzing fMRI data applies some variant of the linear model, which is restrictive in modeling assumptions. To yield more accurate prediction of the time-course behavior of neuronal responses, the semiparametric inference for the underlying hemodynamic response function is developed to identify significantly activated voxels. Under mild regularity conditions, we demonstrate that a class of the proposed semiparametric test statistics, based on the local linear estimation technique, follow $\chi^2$ distributions under null hypotheses for a number of useful hypotheses. Furthermore, the asymptotic power functions of the constructed tests are derived under the fixed and contiguous alternatives. Simulation evaluations and real fMRI data application suggest that the semiparametric inference procedure provides more efficient detection of activated brain areas than the popular imaging analysis tools AFNI and FSL.

**1. Introduction.** Neuroscience is a discipline dedicated to studying the structure, function and pathology of the brain and nervous system, and lies at the forefront of investigation of the brain and mind. Functional magnetic resonance imaging (fMRI) has emerged as a new and exciting noninvasive imaging technique that aims to localize functional brain areas in a living human brain, that is, to detect areas or regions that are responsible for the processing of certain stimuli.

Received May 2007; revised May 2007.

[1]Supported in part by NSF Grants DMS-03-53941, DMS-07-05209 and Wisconsin Alumni Research Foundation.

[2]Supported in part by NSF Grant DMS-03-53941 and Wisconsin Alumni Research Foundation.

Adequate statistical modeling and analysis of the massive spatio-temporal data sets generated by fMRI pose significant challenges to conventional statistical methods. First, a typical fMRI data set for a single scan on a single subject contains a (temporally) highly correlated time series of measurements taken every two seconds or so for about an hour on each of, say, $64 \times 64 \times 30$ voxels (a **voxel** is a volume element in three-dimensional space) throughout the brain. Accordingly, the data sets are so enormous that proper accomodation of both temporal and spatial correlation is needed. Second, models relating fMRI signals to neural changes are complex. The standard tool for analyzing fMRI data is some variant of the linear model, usually fitted separately by least-squares to each voxel [Worsley and Friston (1995)]. After that, tests of significance of the model parameters are performed and colors are drawn on top of significant voxels. This comprises the major procedure of statistical parametric mapping (SPM), popularly used in neuroimage study [Friston et al. (1997)]. Recent reviews of the statistical issues involved in fMRI for brain imaging and the statistical methods for analyzing fMRI data can be found in Lange (1996), Lazar et al. (2001), Fahrmeir and Gössl (2002) and Worsley, et al. (2002), among others.

In this paper, we aim to develop voxelwise semiparametric inference for the underlying hemodynamic response function (HRF), the object of primary interest to neuroscientists. For instance, identifying whether a particular voxel is activated when a subject performs certain motor, sensory or cognitive tasks can be achieved by means of a statistical test of the hypothesis that HRF is zero. In order to generate brain activation maps, statistical inference must be drawn from voxelwise estimates of HRF. We will first develop a semiparametric modeling and estimation approach to obtain statistically more efficient estimates of the underlying HRF associated with fMRI experiments. Compared with the general linear model approach in previous studies, our approach has the advantage that we neither specify any a priori parametric shape for the HRF, nor do we assume any particular form for the temporal drift function. Taking full advantage of these flexibilities will help to reduce the bias due to model misspecification and to enhance the power of detection.

Addressing the issue of semiparametric inference for brain fMRI is a nontrivial task, however. Existing parametric statistical inference procedures for fMRI are not immediately applicable to our approach in which the HRFs are estimated semiparametrically. The work on the generalized likelihood ratio test [Fan, Zhang and Zhang (2001)] sheds light on nonparametric inference, based on function estimation under nonparametric models with independent errors, and, at the same time, is not readily translated into results from other models. Moreover, as emphasized in Section 3, some standard results for semiparametric models are not directly applicable to the context of fMRI data due to the distinctive feature of the Toeplitz design matrix and the complicated dependence structure of the error process. Hence, a rigorous investigation of semiparametric inference applied to the important area of fMRI research is required. This paper fills that gap in the literature. Under mild regularity conditions, we show that a class of the proposed semiparametric test statistics follow $\chi^2$ distributions under null hypotheses for a number of useful hypotheses. To yield improved finite-sample performance of the proposed test statistic, we further explore its bias-corrected version and derive the corresponding asymptotic distribution. Moreover, the asymptotic power functions of the constructed tests are derived under the fixed and contiguous alternatives. These results are not only important for gaining theoretical insight into semiparametric inference applied to a much broader range of scientific problems, but also helpful in offering valuable practical guidance for the implementation of these techniques.

The rest of the paper is arranged as follows. Section 2 reviews statistical models for single-voxel and single-run fMRI. Section 3 describes the semiparametric estimation of the HRF, based on the local linear nonparametric smoothing technique. Section 4 establishes the asymptotic distribution of the proposed test statistics. Section 5 presents simulation evaluations and compares the activated brain regions using the popular imaging analysis tools AFNI (at http://afni.nimh.nih.gov/afni/) [Cox (1996)] and FSL (at http://www.fmrib.ox.ac.uk/fsl/) [Smith et al. (2004) and Woolrich et al. (2001)]. Section 6 applies the semiparametric inference to a real fMRI data set. Technical conditions and detailed proofs are deferred to the Appendix.

**2. Statistical models for single-voxel and single-run fMRI.** We begin with a brief overview of the convolution model popularly used in fMRI study. The BOLD (blood oxygenation level-dependent) signal response to neuronal activity is heavily lagged and damped by the hemodynamic response. Following Ward (2001) and Worsley et al. (2002), a single-voxel fMRI time series $\{s(t_i), y(t_i)\}_{i=1}^n$ for a given scan and a given subject, can be captured by the convolution model,

$$(2.1) \qquad y(t) = s * h(t) + d(t) + \epsilon(t), \qquad t = t_1, \ldots, t_n,$$

where $*$ denotes the convolution operator, $y(t)$ is the measured noisy fMRI signal, $s(t)$ is the external input stimulus at time $t$ [which could be from a design either block- or event-related and where $s(t) = 1$ or $0$ indicates the presence or absence of a stimulus], $h(t)$ is the hemodynamic response function (HRF) at time $t$ after neural activity, $d(t)$ is a slowly drifting baseline of time $t$, and $\epsilon(t)$ is a zero-mean error process, consisting of nonneural noise (due to respiration and blood flow pulsations through the cardiac cycle) and "white noise" (from random/thermal currents in the body and the scanner).

TABLE 1
*HRF, drift and error implemented in AFNI and FSL*

|  | **AFNI (tool** `3dDeconvolve`**)** | **FSL (tool** `FEAT`**)** |
|---|---|---|
| $h(t)$ | finite impulse response filter | difference of two gamma functions, which is the canonical form |
| $d(t)$ | quadratic polynomial | removed in the preprocessing, using high-pass temporal filtering (Gaussian-weighted LSF straight line fitting) |
| $\epsilon(t)$ | i.i.d. | autocorrelation estimated by Tukey tapering of the spectrum of the residuals |

2.1. *Existing methods for modeling HRF, drift and error.* In neuroimaging studies, most existing methods model $h(\cdot)$ as the difference of two gamma functions or a linear combination of gamma functions, a linear combination of a gamma function and its Taylor expansion [Worsley et al. (2002), Lange and Zeger (1997), Josephs and Henson (1999)]. Genovese (2000) constructed $h(\cdot)$ as a "bell" function with cubic splines. As a nuisance component in (2.1), the temporal drift $d(\cdot)$ is usually approximated by a quadratic or higher-order polynomial [Worsley et al. (2002)] or polynomial splines [Genovese (2000)]. Note that restrictive assumptions on the HRF and drift may produce biased estimates of the true hemodynamic responses. Goutte, Nielsen and Hansen (2000) estimated $h(\cdot)$ using smooth FIR filters and reported that some subtle details of the HRF can be revealed by the filters, but not by previous approaches based on gamma functions. The errors $\epsilon(t_i)$ are well known to be temporally autocorrelated. Genovese (2000) assumed independent errors for computational convenience. Other assumptions like the AR($p$) structure, most commonly AR(1), are used in Worsley et al. (2002). As an illustration, Table 1 tabulates the HRF, drift and error implemented in software AFNI and FSL.

**3. Semiparametric estimation of HRF.** Estimating the HRF in (2.1) is a deconvolution problem. Ideally, the HRF is a high-dimensional smooth function and is nonidentically zero if the voxel responds to the stimuli. We will describe a semiparametric method for characterizing properties of the hemodynamic response in the presence of unknown smooth drift. Such characterization is essential for accurate prediction of time-course behavior of neuronal responses.

Typically, the peak value of HRF $h(\cdot)$ is reached after a short delay of the stimulus and drops quickly to zero. A typical example of $h(\cdot)$, given in Glover (1999), is plotted in Figure 1. Clearly, the region $\{t : h(t) \neq 0\}$ is sparse in its temporal domain. Thus, to obtain statistically efficient estimates of the HRF associated with event-related fMRI experiments, the **sparsity** of the HRF needs to be taken into account. We thus suppose that $h(t) = 0$ for $t > t_m$ and focus on estimating the first $m$ values of $h(t_i)$, where $m$ is less than $n$, the length of the fMRI time series. Similarly to the regularization technique discussed in Bickel and Li (2006), such a *qualitative* assumption aims to obtain well-behaved solutions to overparametrized estimation problems and is thus particularly appealing for dimension reduction with high-dimensional problems. The semiparametric modeling and inference in this paper are applicable to all $m < n$. Data-driven selection of $m$ can be made via a change-point approach or other model-selection criteria. To facilitate discussion, we assume that $y(\cdot)$ and $s(\cdot)$ have equal time resolutions of one second. Letting $\mathbf{y} = (y(t_1), \ldots, y(t_n))^T$,

$$\mathbf{S} = \begin{bmatrix} s(0) & 0 & \cdots & 0 \\ s(t_2 - t_1) & s(0) & \cdots & 0 \\ \vdots & \vdots & \ddots & \vdots \\ s(t_m - t_1) & s(t_m - t_2) & \cdots & s(0) \\ \vdots & \vdots & \cdots & \vdots \\ s(t_n - t_1) & s(t_n - t_2) & \cdots & s(t_n - t_m) \end{bmatrix},$$

$\mathbf{h} = (h(t_1), \ldots, h(t_m))^T$, $\mathbf{d} = (d(t_1), \ldots, d(t_n))^T$ and $\boldsymbol{\epsilon} = (\epsilon(t_1), \ldots, \epsilon(t_n))^T$, model (2.1) can be re-expressed as $\mathbf{y} = \mathbf{S}\mathbf{h} + \mathbf{d} + \boldsymbol{\epsilon}$, where $\mathbf{S}$ is a Toeplitz matrix.

In general, for multiple types of stimuli, model (2.1) can be extended to be

$$(3.1) \quad y(t) = s_1 * h_1(t) + \cdots + s_r * h_r(t) + d(t) + \epsilon(t), \qquad t = t_1, \ldots, t_n.$$

Corresponding to the $j$th type of stimulus, denote by $s_j(\cdot)$ the time-varying stimulus function, by $\mathbf{S}_j$ the $n \times m$ Toeplitz design matrix and by $\mathbf{h}_j$ the $m \times 1$ vector of the HRF. Model (3.1) can then be rewritten as

$$(3.2) \qquad \mathbf{y} = \mathbf{S}_1 \mathbf{h}_1 + \cdots + \mathbf{S}_r \mathbf{h}_r + \mathbf{d} + \boldsymbol{\epsilon} \equiv \mathbf{S}\mathbf{h} + \mathbf{d} + \boldsymbol{\epsilon},$$

where $\mathbf{S} = [\mathbf{S}_1, \ldots, \mathbf{S}_r]$ and $\mathbf{h} = [\mathbf{h}_1^T, \ldots, \mathbf{h}_r^T]^T$. To accommodate fMRI data with multiple runs, we only need to supplement the matrix $\mathbf{S}$ by adding the Toeplitz design matrix arising from each run.

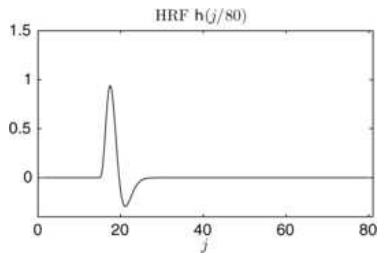

Fig. 1. *An illustrative plot of HRF $h(t_j)$ with $n = 80$.*

Model (3.2) is conceivably a semiparametric regression model, with a vector $\mathbf{h}$ of length $rm$ for *parametric* components and a vector $\mathbf{d}$ of length $n$ for *nonparametric* components. The parametric components (related to the unknown HRF) are of primary interest, whereas the nonparametric components (related to the unknown temporal drift) serve as nuisance effects, and the noise components $\boldsymbol{\epsilon}$ are serially correlated. We wish to emphasize that due to the special structure of the design matrix $\mathbf{S}$ associated with fMRI design, some commonly used assumptions, such as independence between rows of a design matrix, fail to hold. In addition, the unobservable true correlation structure of $\boldsymbol{\epsilon}$ is often complicated. Thus some standard results for semiparametric models are not directly applicable to the current fMRI data.

We now describe the semiparametric estimation of both the HRF and the nonparametric drift function in (3.2). Let $S_d$ be an $n \times n$ local linear smoothing matrix associated with the design points $\{t_1, \ldots, t_n\}$, with the $(i, j)$th entry equal to

$$(3.3) \quad S_d(i,j) = (1,0)\{\mathbf{X}(t_i)^T \mathbf{W}(t_i)\mathbf{X}(t_i)\}^{-1}(1, t_j - t_i)^T K((t_j - t_i)/b)/b,$$

where $K$ is a kernel function, $b > 0$ is a bandwidth parameter,

$$\mathbf{X}(t) = \begin{bmatrix} 1 & t_1 - t \\ \vdots & \vdots \\ 1 & t_n - t \end{bmatrix}$$

and

$$\mathbf{W}(t) = \operatorname{diag}\{K((t_1 - t)/b)/b, \ldots, K((t_n - t)/b)/b\};$$

see Fan and Gijbels (1996), which provides a comprehensive account of the local linear and local polynomial regression techniques. (For expositional simplicity, this paper is confined to the local linear method.) Note that the matrix $S_d$ carries information about the design points, kernel $K$ and bandwidth $b$, but does not rely on the configuration of the response variables. We refer to Section 2.3 of Zhang (2003) for further discussion of finite-sample and asymptotic properties of the smoothing matrix. Notice that smoothing the entries of $\mathbf{y}$ via the local linear method is equivalent to applying $S_d$ to $\mathbf{y}$. We observe from (3.2) that

$$(3.4) \quad \widetilde{\mathbf{y}} = \widetilde{\mathbf{S}}\mathbf{h} + \widetilde{\mathbf{d}} + \widetilde{\boldsymbol{\epsilon}},$$

where $\widetilde{\mathbf{y}} = (\mathbf{I} - S_d)\mathbf{y}$, $\widetilde{\mathbf{S}} = (\mathbf{I} - S_d)\mathbf{S}$, $\widetilde{\mathbf{d}} = (\mathbf{I} - S_d)\mathbf{d}$, $\widetilde{\boldsymbol{\epsilon}} = (\mathbf{I} - S_d)\boldsymbol{\epsilon}$ and $\mathbf{I}$ denotes an identity matrix. Ignoring $\widetilde{\mathbf{d}}$, model (3.4) can be regarded as a general linear model. Denote by $R$ the true correlation matrix of $\boldsymbol{\epsilon}$, namely, $\operatorname{cov}(\boldsymbol{\epsilon}, \boldsymbol{\epsilon}) = \sigma^2 R$, with variance $\sigma^2$. Let $\widehat{R}$ be an estimate of $R$. By the weighted least-squares method, an estimate of $\mathbf{h}$ is produced by

$$(3.5) \quad \widehat{\mathbf{h}} = (\widetilde{\mathbf{S}}^T \widehat{R}^{-1} \widetilde{\mathbf{S}})^{-1} \widetilde{\mathbf{S}}^T \widehat{R}^{-1} \widetilde{\mathbf{y}},$$

which, in turn, supplies estimates of the drift components formed by

$$\widehat{\mathbf{d}} = S_d(\mathbf{y} - \mathbf{S}\widehat{\mathbf{h}}).$$

**4. Semiparametric hypothesis test for HRF.** Identification of a particular brain region with a specific function has become a central theme in neuroscience. In this section, we consider constructing test statistics to test whether a particular voxel is activated by the stimuli and whether HRFs activated by different types of stimuli really differ. They correspond to testing the hypotheses $H_0 : \mathbf{h} = \mathbf{0}$ versus $H_1 : \mathbf{h} \neq \mathbf{0}$ and $H_0 : \mathbf{h}_{j_1} = \mathbf{h}_{j_2}$ versus $H_1 : \mathbf{h}_{j_1} \neq \mathbf{h}_{j_2}$, where $j_1 \neq j_2$. Under the semiparametric model (3.2), all of these testing problems can be formulated in a more general form,

(4.1) $\qquad\qquad H_0 : A\mathbf{h} = \mathbf{0} \quad \text{versus} \quad H_1 : A\mathbf{h} \neq \mathbf{0},$

where $A$ is a full row rank matrix with $\text{rank}(A) = \mathsf{k}$.

An earlier work on developing pseudo-$F$-type test statistics was empirically studied in Lu (2006) and Zhang et al. (2006). There, it was observed from QQ plots in simulation studies that under the null hypothesis, empirical quantiles of the $F$-type test statistics (in the restrictive case where the true $R$ is known and a single type of stimulus in the fMRI experiment is presented) could be approximated by quantiles of the $F$-distribution. No asymptotic exploration of properties of the $F$-type test statistics was conducted.

4.1. *Asymptotic null distributions.* Motivated by the parametric $F$-statistic in linear regression models and the justification of power comparison [Zhang and Dette (2004)] between nonparametric tests for regression curves based on kernel smoothing techniques, we first examine the following semiparametric test statistic, represented by

$$\mathbb{K} = \frac{(A\widehat{\mathbf{h}})^T \{A(\widetilde{\mathbf{S}}^T \widehat{R}^{-1} \widetilde{\mathbf{S}})^{-1} A^T\}^{-1}(A\widehat{\mathbf{h}})}{\widehat{\mathbf{r}}^T \widehat{R}^{-1} \widehat{\mathbf{r}}/(n - rm)},$$

where $\widehat{\mathbf{r}} = \widetilde{\mathbf{y}} - \widetilde{\mathbf{S}}\widehat{\mathbf{h}}$. Theorem 4.1 below establishes the asymptotic null distribution of $\mathbb{K}$.

THEOREM 4.1. *Assume Condition* A *in the Appendix. Then, under $H_0$ in (4.1), where $A$ is a $\mathsf{k} \times rm$ matrix with $\text{rank}(A) = \mathsf{k}$, it follows that $\mathbb{K} \xrightarrow{\mathcal{L}} \chi^2_{\mathsf{k}}$, where $\xrightarrow{\mathcal{L}}$ denotes convergence in distribution.*

Our simulation evaluation in Section 5 demonstrates that the finite sampling distribution of $\mathbb{K}$ is reasonably well approximated by its asymptotic $\chi^2$ distribution, whereas when the noise level decreases, the approximation may become less accurate; see Figure 2 (right panel). Technically, as manifested in the proof of Theorem 4.1, the asymptotic $\chi^2$ distribution of $\mathbb{K}$ follows from the asymptotic normality of $\widehat{\mathbf{h}}$ shown in Lemma A.7, which relies on the fact that a term $J_1$ (associated with the drift vector $\mathbf{d}$) is stochastically dominated by a term $J_2$ (associated with the error vector $\boldsymbol{\epsilon}$). Practically, in finite-sample situations, low noise levels do not necessarily guarantee that $J_1$ is stochastically negligible compared with $J_2$. Consequently, the finite sampling distributions of $\widehat{\mathbf{h}}$ and $\mathbb{K}$ may appear biased toward the normal and $\chi^2$ distribution, respectively. In these situations, we adopt the bias-corrected version of $\mathbb{K}$, defined as

$$\mathbb{K}_{\mathrm{bc}} = \frac{(A\widehat{\mathbf{h}}_{\mathrm{bc}})^T \{A(\widetilde{\mathbf{S}}^T \widehat{R}^{-1} \widetilde{\mathbf{S}})^{-1} A^T\}^{-1} (A\widehat{\mathbf{h}}_{\mathrm{bc}})}{\widehat{\mathbf{r}}_{\mathrm{bc}}^T \widehat{R}^{-1} \widehat{\mathbf{r}}_{\mathrm{bc}}/(n-rm)},$$

where $\widehat{\mathbf{h}}_{\mathrm{bc}} = \widehat{\mathbf{h}} - (\widetilde{\mathbf{S}}^T \widehat{R}^{-1} \widetilde{\mathbf{S}})^{-1} \widetilde{\mathbf{S}}^T \widehat{R}^{-1} \widetilde{\widehat{\mathbf{d}}}$, $\widehat{\mathbf{r}}_{\mathrm{bc}} = \widehat{\mathbf{r}} - \widetilde{\widehat{\mathbf{d}}}$, $\widehat{\mathbf{d}} = S_d(\mathbf{y} - \mathbf{S}\widehat{\mathbf{h}})$ and $\widetilde{\widehat{\mathbf{d}}} = (\mathbf{I} - S_d)\widehat{\mathbf{d}}$. Note that as the sequence length $n$ grows, $\widetilde{\widehat{\mathbf{d}}}$ is negligible, but practically adjusts for the bias caused by $J_1$ due to the ignorance of $\widetilde{\mathbf{d}}$ in (3.4). Theorem 4.2 below reveals that $\mathbb{K}_{\mathrm{bc}}$ and $\mathbb{K}$ have the same asymptotic null distributions.

THEOREM 4.2. *Assume Condition* A *in the Appendix. Then, under* $H_0$ *in (4.1), where $A$ is a* $\mathsf{k} \times rm$ *matrix with* $\mathrm{rank}(A) = \mathsf{k}$*, it follows that* $\mathbb{K}_{\mathrm{bc}} \overset{\mathcal{L}}{\to} \chi^2_{\mathsf{k}}$.

We now make some remarks concerning the derivations of Theorems 4.1–4.2. First, it is tempting to try to show that $n^{-1}\widetilde{\mathbf{S}}^T \widehat{R}^{-1} \widetilde{\mathbf{S}} \overset{\mathcal{P}}{\to} \mathbf{M}$ for some positive definite matrix $\mathbf{M}$, where $\overset{\mathcal{P}}{\to}$ denotes converges in probability. Nonetheless, for fMRI data, since the $n \times n$ correlation matrix $R$ of $\boldsymbol{\epsilon}$ is generally far more complicated than the diagonal matrix of independent errors, deriving an explicit form for $\mathbf{M}$ is nearly intractable. To overcome this technical difficulty, we have demonstrated that it suffices to verify that $R$ satisfies

$$\mathrm{var}\{n^{-1}\boldsymbol{\xi}_{j_1,\ell_1}^T (\mathbf{I} - S_d)^T R^{-1} (\mathbf{I} - S_d) \boldsymbol{\xi}_{j_2,\ell_2}\} \to 0$$

for all $j_1, j_2 = 1, \ldots, r$ and all $\ell_1, \ell_2 = 1, \ldots, m$, where $\boldsymbol{\xi}_{j,\ell}$ is the $\ell$th column vector of $\mathbf{S}_j$ and $\widehat{R}$ fulfills Condition A8 in the Appendix,

$$E(\|\widehat{R}^{-1} - R^{-1}\|_\infty^2) = o(1),$$

$\|B\|_\infty = \max_{1 \le i \le n} \sum_{j=1}^n |B(i,j)|$ denoting the $\infty$-norm of an $n \times n$ matrix $B$; see Lemma A.6 and Corollary A.2. Thus, the explicit form of $\mathbf{M}$ is not needed in deriving the asymptotic null distributions of $\mathbb{K}$ and $\mathbb{K}_{\mathrm{bc}}$. Second, Condition A8, together with $\|B\|_2 \le \{\|B\|_1 \|B\|_\infty\}^{1/2}$ [Golub and Van Loan (1996)] and the symmetry of $\widehat{R}$ and $R$, guarantees that $\|\widehat{R}^{-1} - R^{-1}\|_2 = o_P(1)$, which is typically interpreted as the "consistency" of large covariance matrix estimators [Bickel and Levina (2008)].

REMARK 1. In real-world applications, fMRI sequence lengths are not very long. For instance, $n$ is 185 for each run in the real fMRI data set described in Section 6. This indicates that the "mixing assumptions," commonly made in the asymptotic studies of nonlinear time series [Bosq (1998), Fan and Yao (2003)], may not hold for fMRI data. Therefore, the sampling properties of $\mathbb{K}$ and $\mathbb{K}_{\mathrm{bc}}$ are studied using the more realistic error assumption A3 of Condition A in the Appendix, which could possibly be weakened.

REMARK 2. Throughout the numerical work in this paper, parametric estimation of the error covariance matrix adopts a computationally fast and effective scheme developed in Zhang et al. (2006), which assumes $g = 2$ in Condition A3 of the Appendix. This regularized estimator is constructed as follows: obtain the transformed data $e(t_i)$ by applying the second-order difference to $y(t_i) - \sum_{j=1}^{r} s_j * h_j(t_i)$; calculate autocovariances $\{\gamma_e(j)\}_{j=0}^{g}$ of $e(t_i)$, which form a linear system for autocovariances $\{\gamma(j)\}_{j=0}^{g}$ of $\epsilon(t_i)$; substitute for $\{\gamma_e(j)\}$ their empirical moment estimates and solve $\{\gamma(j)\}$; acquire an estimate $\widehat{R}$ of $R$ using Condition A3. Moreover, since an fMRI data set contains time-course measurements over voxels, the number of which is typically of the order of $10^4$–$10^5$, the conventional false discovery rate (FDR) approach [Benjamini and Hochberg (1995), Storey (2002)] can be adopted to account for the multiple comparisons problem. Other useful and elaborate procedures for covariance matrix estimation and multiple comparison may also be employed. Particularly, Zhang et al. (2006) presented numerical evidence that the existing FDR approach tends to find activation in tiny scattered regions of the brain which are more likely to be false discoveries, and carefully devised a new FDR approach which gains efficiency over the existing FDR approach.

4.2. *Asymptotic power functions.* To appreciate the discriminating power of the proposed tests in assessing the significance of activated areas, the asymptotic power is analyzed. Theorem 4.3 demonstrates that both $\mathbb{K}$ and $\mathbb{K}_{\mathrm{bc}}$ are consistent against all fixed deviations from the null model.

THEOREM 4.3. *Assume Condition* A *in the Appendix and* $n^{-1}\widetilde{\mathbf{S}}^T R^{-1} \widetilde{\mathbf{S}} \xrightarrow{\mathcal{P}} \mathbf{M}$, *where* $\mathbf{M}$ *is positive definite. Then, under the fixed alternative* $H_1$ *in* (4.1),

$$n^{-1}\mathbb{K} \xrightarrow{\mathcal{P}} (A\mathbf{h})^T (A\mathbf{M}^{-1}A^T)^{-1} A\mathbf{h}/\sigma^2 > 0,$$
$$n^{-1}\mathbb{K}_{\mathrm{bc}} \xrightarrow{\mathcal{P}} (A\mathbf{h})^T (A\mathbf{M}^{-1}A^T)^{-1} A\mathbf{h}/\sigma^2 > 0.$$

The results in Theorem 4.3 indicate that under the fixed alternative $H_1$,

$$\mathbb{K} \xrightarrow{\mathcal{P}} +\infty \quad \text{and} \quad \mathbb{K}_{\text{bc}} \xrightarrow{\mathcal{P}} +\infty,$$

at the common rate $n$. Hence, the test statistics $\mathbb{K}$ and $\mathbb{K}_{\text{bc}}$ have power functions tending to one against fixed alternatives.

Consider a sequence of local alternatives, defined by

$$(4.2) \qquad H_{1n} : A\mathbf{h} = \delta_n \mathbf{c},$$

where $\delta_n = n^{-1/2}$ and $\mathbf{c} = (c_1, \ldots, c_{\mathsf{k}})^T \neq \mathbf{0}$. Theorem 4.4 explores the asymptotic distributions of $\mathbb{K}$ and $\mathbb{K}_{\text{bc}}$ under the local alternatives $H_{1n}$.

THEOREM 4.4. *Assume Condition* A *in the Appendix and* $n^{-1}\widetilde{\mathbf{S}}^T R^{-1} \times \widetilde{\mathbf{S}} \xrightarrow{\mathcal{P}} \mathbf{M}$, *where* $\mathbf{M}$ *is positive definite. Then, under the local alternative* $H_{1n}$ *in* (4.2), $\mathbb{K} \xrightarrow{\mathcal{L}} \chi^2_{\mathsf{k}}(\tau^2)$ *and* $\mathbb{K}_{\text{bc}} \xrightarrow{\mathcal{L}} \chi^2_{\mathsf{k}}(\tau^2)$, *with noncentrality parameter* $\tau^2 = \mathbf{c}^T (A\mathbf{M}^{-1}A^T)^{-1}\mathbf{c}/\sigma^2$.

The results in Theorem 4.4 indicate that the tests have nontrivial local power detecting local alternatives approaching the null at the rate $n^{-1/2}$. A simple calculation shows that the asymptotic power of the tests against local misspecification (4.2) equals

$$\int_{\chi^2_{\mathsf{k};1-\alpha}}^{\infty} \frac{\exp\{-(x+\tau^2)/2\}}{2^{\mathsf{k}/2}} \sum_{j=0}^{\infty} \frac{x^{\mathsf{k}/2+j-1}\tau^{2j}}{\Gamma(\mathsf{k}/2+j)2^{2j}j!} \, dx,$$

where $\chi^2_{\mathsf{k};1-\alpha}$ is the $1-\alpha$ quantile of the $\chi^2_{\mathsf{k}}$ distribution and $\Gamma(\cdot)$ denotes the gamma function.

**5. Simulation study.** Throughout the numerical work, we use the Epanechnikov kernel function [Silverman (1986)] supported on $[-1, 1]$. A complete copy of MATLAB codes is available on request.

5.1. *Hypothesis test of HRF at a single voxel.* As an illustration, the hypothesis testing for $H_0 : \mathbf{h} = \mathbf{0}$ versus $H_1 : \mathbf{h} \neq \mathbf{0}$ is undertaken. This is used to test whether the brain activity in a voxel is triggered or not. To check the agreement between the $\chi^2$ distribution with finite sampling distributions of $\mathbb{K}$ and $\mathbb{K}_{\text{bc}}$ under $H_0$, the fMRI data are simulated as follows. We simulate an fMRI experiment with a single run and a single type of stimulus, where $n = 400$ and 500 realizations are conducted. **(I)** The time-varying stimuli are generated from independent Bernoulli trials such that $P\{s(t_i) = 1\} = 0.5$. **(II)** The HRF is $h(t_i) = 0$, $i = 1, \ldots, 18$ (so that $m = 18$). **(III)** The drift function is $d(t_i) = 10\sin\{\pi(t_i - 0.21)\}$, $i = 1, \ldots, n$. **(IV)** The noise process $\epsilon$ is the sum of independent noise processes $\epsilon_1$ and $\epsilon_2$ (see Purdon et al. (2001));

$\{\epsilon_1(t_i)\}$ are i.i.d. normal with mean zero and variances $0.5216^2$, $0.3689^2$, $0.2608^2$ and $0.1844^2$, respectively; $\epsilon_2$ is AR(1), that is, $\epsilon_2(t_i) = \rho\epsilon_2(t_{i-1}) + z(t_i)$ with $\rho = 0.638$ and the $z(t_i)$ follow the normal distribution with mean zero and variances $0.5216^2$, $0.3689^2$, $0.2608^2$ and $0.1844^2$, respectively. These choices give a noise lag-one autocorrelation equal to 0.4 and signal-to-noise-ratios (SNRs) of about 1, 2, 4 and 8, where SNR = variance($\mathbf{Sh}$)/variance($\boldsymbol{\epsilon}$).

The QQ plots of the (1st to 99th) percentiles of $\mathbb{K}$ and $\mathbb{K}_{\mathrm{bc}}$ against those of the $\chi^2_m$ distribution are displayed in Figure 2. In the top panel, $\mathbb{K}$ and $\mathbb{K}_{\mathrm{bc}}$ use the true covariance matrix and fix the smoothing parameters at their theoretically optimal values (minimizing the mean squared errors of estimators) for estimating the HRF and drift in each simulation. For the sake of clarity, only the cases of SNR equal to 1 and 8 are presented; the former is the "large noise level" case, whereas the latter is the "small noise level" case. In either case, we observe that the finite sampling distributions of $\mathbb{K}$ and $\mathbb{K}_{\mathrm{bc}}$, at the realistic sample size 400, agree reasonably well with the $\chi^2$ distribution. The QQ plots also lend support to the possibility that $\mathbb{K}_{\mathrm{bc}}$ is better than or at least as good as, the bias-uncorrected counterpart $\mathbb{K}$.

For a more realistic comparison, $\mathbb{K}$ and $\mathbb{K}_{\mathrm{bc}}$ in the bottom panel of Figure 2 use the estimated covariance matrices and data-driven smoothing parameters. The results are similar in spirit to the ones in the top panel and continue to support the bias correction procedure.

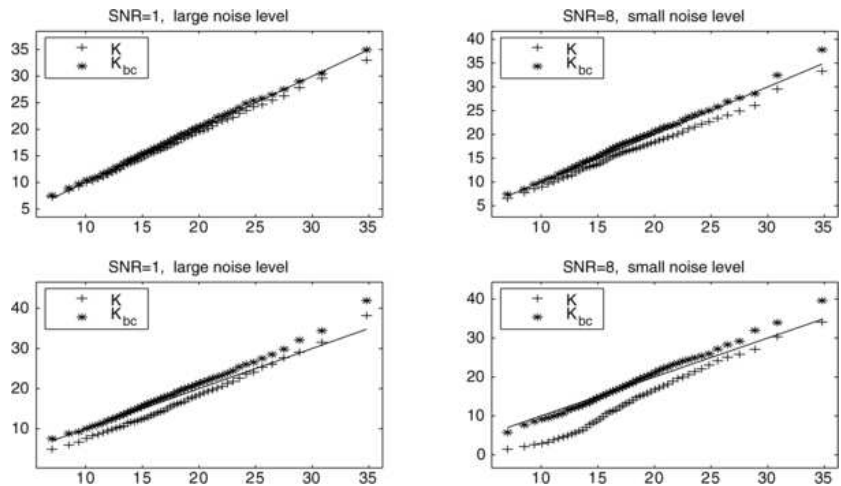

FIG. 2. *Empirical quantiles (on the y-axis) of test statistics $\mathbb{K}$ and $\mathbb{K}_{\mathrm{bc}}$ (where the top panel uses the true $R$ and the optimal smoothing parameters, and the bottom panel uses the estimated $\widehat{R}$ and data-driven smoothing parameters) versus quantiles (on the x-axis) of $\chi^2_m$ distribution. Solid line: the 45 degree reference line.*

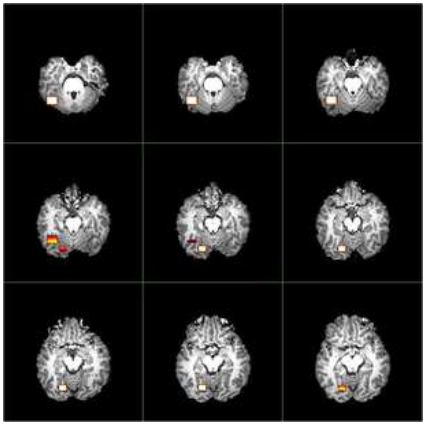

Fig. 3. *True activated brain regions (denoted by hot color) for the simulated fMRI data set.*

5.2. *Detection of activated brain regions.* We simulate a whole brain fMRI data set, with aim of mimicking true brain activity to the maximum extent feasible. The experiment design, timings and size are exactly the same as those of the real fMRI data set in Section 6. An HRF profile is extracted from a voxel which shows the strongest responses in the real data set. For each voxel, the simulated drift is obtained from an adequate smoothing of the time series for the corresponding voxel of the real data set. The simulated noise variance profile is determined from a variance map, which is made by a $5 \times 5 \times 5$ spatial median smoothing on median values of squared residuals of the real time series, subtracting the simulated drift profile as mentioned before. The noise process $\epsilon(t)$ is generated in a fashion similar to that of Section 5.1. Specifically, the variances of $\epsilon_1(t)$ and $z(t)$ are chosen to be equal such that var$\{\epsilon(t)\}$ is one-fifth of the variance map. The HRF profiles, in accordance with the stimuli in the experiment, is added to two regions which are postulated to be truly active. In these two zones, the HRFs have been rescaled to about 17% and 12% of the amplitude of the original HRF profiles. The purpose of rescaling the HRFs and noise variance is to amplify the drift effect and weaken the HRF response so that the estimation of the HRF is more challenging. Figure 3 shows nine different slices which highlight the two activated brain regions. Note that throughout the paper, we apply the same registration transform from the real brain data to the T1 high-resolution image of the subject's brain.

The gain in efficiency achieved by the semiparametric inference procedure is illustrated by comparing the activated brain regions identified by our approach with those identified via methods offered by the popular software AFNI and FSL. The conventional FDR approach is performed at the FDR level 0.05. Inspection of Figure 4 reveals that $\mathbb{K}$ and $\mathbb{K}_{bc}$ are capable of locating both active regions. In contrast, both AFNI and FSL fail to locate an activated brain area, and the other region, although correctly detected, has appreciably reduced size relative to the actual size. This detection bias suggests that the stringent modeling assumptions in Table 1 should be relaxed to ameliorate the effects of misspecification. Furthermore, as evidenced in Figure 5, all four methods, when applying the new FDR approach in Zhang et al. (2006), achieve more accurate detection than their counterparts in Figure 4, with $\mathbb{K}$ and $\mathbb{K}_{bc}$ continuing to outperform AFNI and FSL. Therefore, for applications to the real fMRI data set in Section 6, we will only employ the new FDR approach in Zhang et al. (2006).

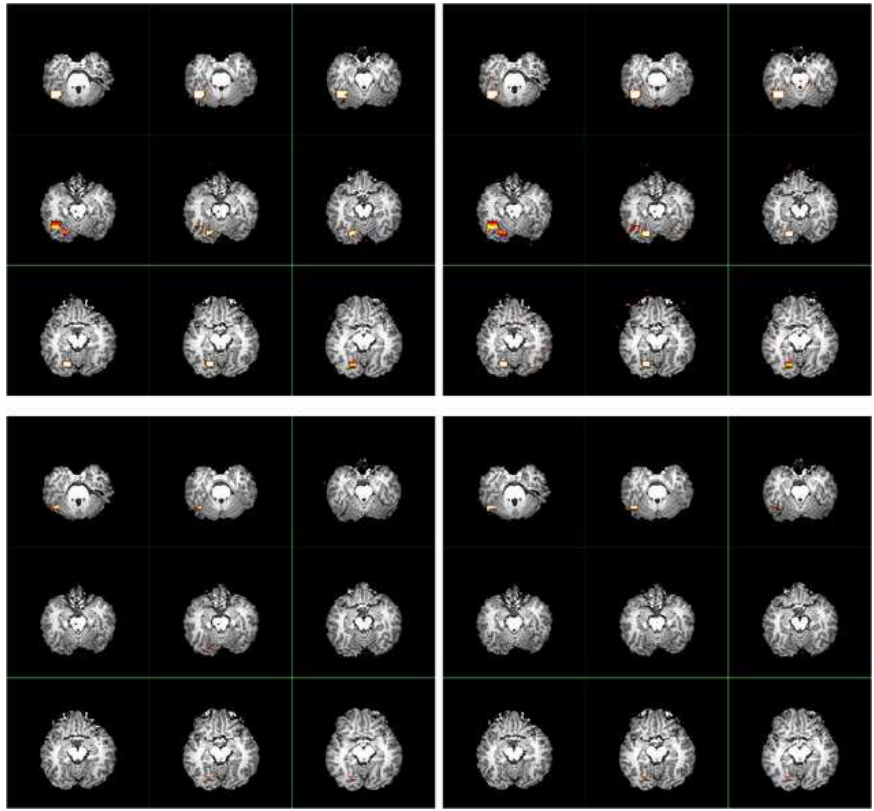

Fig. 4. *Comparison of activated brain regions discovered for the simulated fMRI data set. Top panel:* $\mathbb{K}$ (*on the left*) *and* $\mathbb{K}_{bc}$ (*on the right*). *Bottom panel: AFNI* (*on the left*) *and FSL* (*on the right*). *The conventional FDR approach is used. The FDR level is* 0.05.

**6. Real data analysis.** In an emotional control study, subjects saw a series of negative or positive emotional images and were asked to either suppress or enhance their emotional responses to the image, or to simply attend to the image. Therefore, there were six types of trial (i.e., six types of stimuli): negative-enhance (neg-enh), negative-attend (neg-att), negative-suppress (neg-sup), positive-enhance (pos-enh), positive-attend (pos-att) and positive-suppress (pos-sup). The sequence of trials was randomized. The time between successive trials also varied. There were 24 trials each of neg-enh, neg-sup, pos-enh, and pos-sup; there were 11 trials each of neg-att and pos-att.

The size of the whole brain data set is $64 \times 64 \times 30$. At each voxel, the time series has six runs, each containing 185 observations with a time resolution of two seconds, thus TR = 2 seconds and the total length is 1110. In contrast, the length of stimuli is 2220; the timing of the stimuli has a time resolution of one second and thus each HRF output will also be sampled at one second.

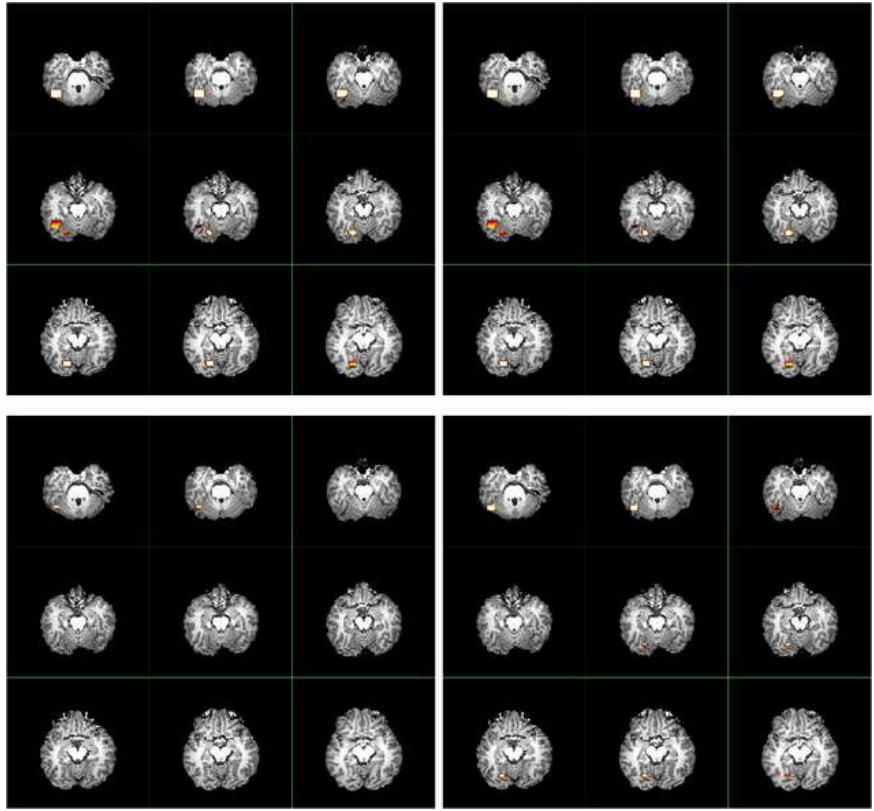

FIG. 5. *Comparison of activated brain regions discovered for the simulated fMRI dataset. Top panel: $\mathbb{K}$ (on the left) and $\mathbb{K}_{bc}$ (on the right). Bottom panel: AFNI (on the left) and FSL (on the right). The new FDR approach in Zhang et al. (2006) is used. The FDR level is* 0.05.

SEMI-PARAMETRIC DETECTION FOR FMRI 15

Hence, the odd rows of the design matrix **S** in (3.2) suffice for analysis. The study aims to estimate the BOLD response to each of the trial types for 1–18 seconds following the image onset. We analyze the fMRI data set containing one subject. The length of the estimated HRF is set equal to 18.

A comparison of the activated brain regions detected by $\mathbb{K}$, $\mathbb{K}_{\mathrm{bc}}$, AFNI and FSL is illustrated in Figure 6. Again, the HRF in FSL is specified as the difference of two gamma functions and the drift term in AFNI is specified as a quadratic polynomial. We use FDR at level 0.001 to carry out the multiple comparisons. This level is set to avoid excessive discoveries, most of which are thought to be false. Our detected regions are closer to those obtained by AFNI, but our methods find activation in much more clustered regions of the brain. For example, our results do not have the holes seen in the detected regions on the first slice of AFNI and FSL. AFNI gives more tiny scattered findings, which are more likely to be false discoveries. FSL detects very scattered regions which are difficult to interpret. In addition, the volumes of the regions detected by FSL are substantially smaller than those of AFNI and our methods.

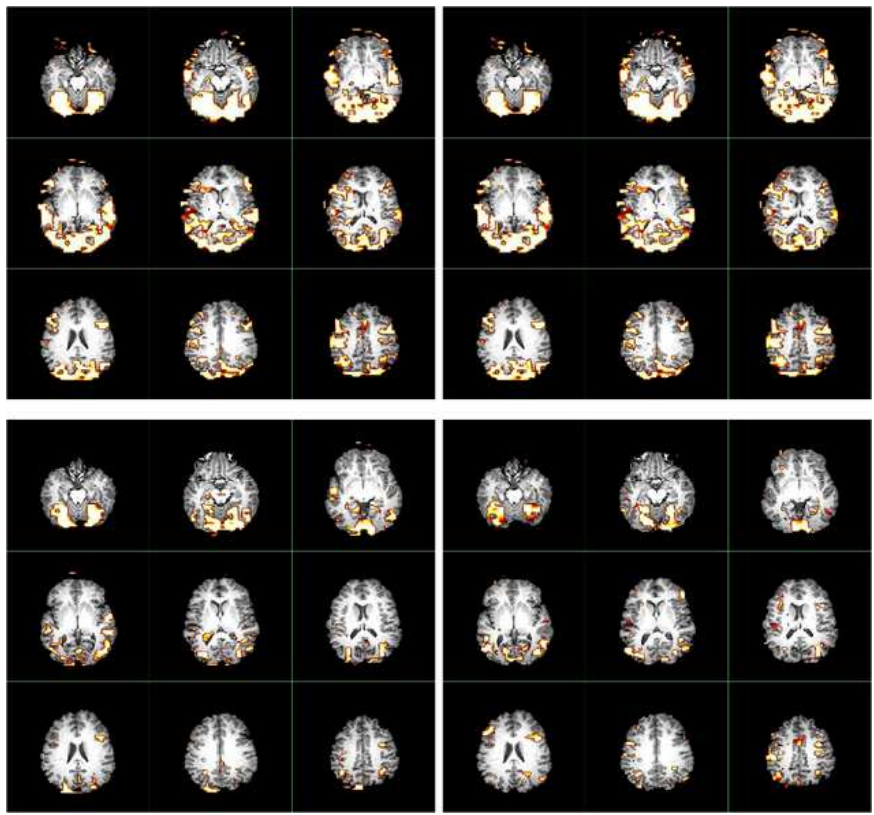

FIG. 6. *Comparison of activated brain regions discovered for the real fMRI data set. Top panel:* $\mathbb{K}$ (*on the left*) *and* $\mathbb{K}_{\mathrm{bc}}$ (*on the right*). *Bottom panel: AFNI* (*on the left*) *and FSL* (*on the right*). *The new FDR approach in Zhang et al.* (*2006*) *is used. The FDR level is* 0.001.

## APPENDIX: PROOFS OF MAIN RESULTS

We first impose some technical assumptions, which are not the weakest possible. Throughout the proof, $C$ is used as a generic finite constant.

CONDITION A.

A1. The drift function $d(t)$ has a bounded continuous second derivative.

A2. The kernel $K$ is a symmetric probability density function with compact support, say $[-L, L]$, is Lipschitz continuous and such that $\sup_t K(t) \leq C$ for some constant $C \in (0, \infty)$.

A3. Assume that $\{\epsilon(t_i)\}$ is a stationary $g$-dependent sequence with $E\{\epsilon(t_1)\} = 0$, $E\{\epsilon^2 \in (0,\infty)(t_1)\} = \sigma^2$ and $E\{\epsilon^4(t_1)\} < \infty$. The eigenvalues of $R$, the true correlation matrix of $\epsilon$, are uniformly bounded away from zero and infinity. Furthermore, $E\{\epsilon(t_i)\epsilon(t_j)|\widehat{R}\} = E\{\epsilon(t_i)\epsilon(t_j)\}$.

A4. In model (3.1), $\{s_j(\cdot)\}$, $j = 1, \ldots, r$, are independent of $\{\epsilon(\cdot)\}$. For the RPER design, $s_j(t)$ is stationary and $P\{s_j(t) = 1\} = p_j \in (0,1)$, $j = 1, \ldots, r$, and $\sum_{j=1}^r p_j < 1$. Assume that $s_{j_1}(t_u)$ and $s_{j_2}(t_v)$ are independent at any $t_u \neq t_v$. For any $u, v = 1, \ldots, n$, $E\{s_j(t_u)s_j(t_v)|\widehat{R}\} = E\{s_j(t_u)s_j(t_v)\}$.

A5. $n \to \infty$, $b \to 0$ and $nb \to \infty$.

A6. $t_i = i/n$, $i = 1, \ldots, n$.

A7. $\operatorname{cov}(\mathbf{S}^T, \mathbf{S}^T) > \mathbf{0}$.

A8. $E(\|\widehat{R}^{-1} - R^{-1}\|_\infty^2) = o(1)$.

We next introduce some necessary notation and definitions.

NOTATION. For the kernel $K$ and bandwidth $b > 0$, define $K_b(t) = K(t/b)/b$. Denote by $\mathbf{e}_j$ the $j$th column of an identity matrix. Define vectors $\mathbf{1} = (1, \ldots, 1)^T$ and $\mathbf{0} = (0, \ldots, 0)^T$. Define a matrix $H$ with entries $H(i,j) = n^{-1} K_b(t_j - t_i)$, $1 \leq i, j \leq n$. Define $V = R^{-1}$ and let $\rho(l)$ denote the noise autocorrelation coefficient. Denote by $\boldsymbol{\xi}_{j,\ell}$ the $\ell$th column vector of $\mathbf{S}_j$, that is, $\boldsymbol{\xi}_{j,\ell} = \mathbf{S}_j \mathbf{e}_\ell$. Throughout the proof, $\|\cdot\|$ refers to the $L_2$-norm unless otherwise stated.

DEFINITION A.1. An $n \times n$ matrix $B$ is called "row absolute value uniformly summable" (RAVUS) if there exists $C > 0$ such that

$$\sup_{n \geq 1} \sup_{1 \leq j \leq n} \sum_{i=1}^n |B(i,j)| \leq C.$$

Likewise, $B$ is called "column absolute value uniformly summable" (CAVUS) if there exists $C > 0$ such that $\sup_{n \geq 1} \sup_{1 \leq i \leq n} \sum_{j=1}^n |B(i,j)| \leq C$. Moreover, if a matrix is both RAVUS and CAVUS, it is called "absolute value uniformly summable" (AVUS).

Before proving the main results of the paper, we need Lemmas A.1–A.8.

LEMMA A.1. *If both matrices $B_1 \in \mathbb{R}^{n \times n}$ and $B_2 \in \mathbb{R}^{n \times n}$ are AVUS, where AVUS is defined in Definition A.1 above, then $B_1 B_2$ is AVUS.*

PROOF. By the definition, there exists $C > 0$ such that

$$\sup_{n \geq 1} \sup_{1 \leq j \leq n} \sum_{i=1}^n |B_1(i,j)| \leq C, \qquad \sup_{n \geq 1} \sup_{1 \leq j \leq n} \sum_{i=1}^n |B_2(i,j)| \leq C.$$

We observe that

$$\sup_{n \geq 1} \sup_{1 \leq j \leq n} \sum_{i=1}^n |(B_1 B_2)(i,j)| \leq \sup_{n \geq 1} \sup_{1 \leq j \leq n} \sum_{i=1}^n \sum_{l=1}^n |B_1(i,l) B_2(l,j)|$$

$$= \sup_{n \geq 1} \sup_{1 \leq j \leq n} \sum_{l=1}^n |B_2(l,j)| \sum_{i=1}^n |B_1(i,l)|$$

$$\leq C \sup_{n \geq 1} \sup_{1 \leq j \leq n} \sum_{l=1}^n |B_2(l,j)| \leq C^2.$$

Thus, $B_1 B_2$ is RAVUS and, by similar reasoning, is CAVUS. Hence, $B_1 B_2$ is AVUS. □

LEMMA A.2. *Assume Condition* A3. *There then exist constants $C \in (0, \infty)$ and $\lambda \in (0, 1)$ such that $|V(i,j)| \leq C \lambda^{|i-j|}$ for all $1 \leq i, j \leq n$ and $n \geq 1$.*

PROOF. Under Condition A3, $R$ is positive definite, centered and $2g$-banded. Let $a_n$ and $b_n$ be the minimum and maximum eigenvalues of $R$, respectively, and set $r_n = b_n / a_n$. Applying Proposition 2.2 of Demko, Moss and Smith (1984) gives that $|V(i,j)| \leq C_n \lambda_n^{|i-j|}$, where $C_n = \max\{a_n^{-1}, (1 + r_n^{1/2})^2 / (2 a_n r_n)\}$ and $\lambda_n = \{(r_n^{1/2} - 1)/(r_n^{1/2} + 1)\}^{1/g}$. From Condition A3, $a_n$ and $b_n$ are bounded away from both zero and infinity; it follows that $r_n$ is bounded away from both zero and infinity. Thus, there exist $C \in (0, \infty)$ and $\lambda \in (0, 1)$ such that $C_n < C$ and $\lambda_n < \lambda$. Hence, $|V(i,j)| \leq C_n \lambda_n^{|i-j|} \leq C \lambda^{|i-j|}$. □

LEMMA A.3. *Assume Conditions* A2 *and* A5.

1. *Let $f$ be Lipschitz continuous and bounded on an interval $[d_1, d_2]$, where $d_1 < d_2$. Let $u_j = d_1 + (d_2 - d_1)j/n$, $j = 1, \ldots, n$. Then, uniformly in $\tau$,*

$$\text{(A.1)} \quad \frac{1}{n}\sum_{j=1}^{n} K_b(u_j - \tau)f(u_j) = \frac{1}{d_2 - d_1}\int_{d_1}^{d_2} K_b(u - \tau)f(u)\,du + O\left(\frac{1}{nb}\right).$$

2. *Let $\{\epsilon_j\}_{j=1}^{\infty}$ be a sequence of $g$-dependent and identically distributed random variables. Assume $E(\epsilon_1^2) < \infty$. Then, for $t_j = j/n$, $j = 1, \ldots, n$,*

$$\text{(A.2)} \quad \sup_{t \in [bL, 1-bL]} E\left[\left\{\frac{1}{n}\sum_{j=1}^{n}\epsilon_j K_b(t_j - t) - E(\epsilon_1)\right\}^2\right] = O\left(\frac{1}{nb}\right).$$

PROOF. We first show (A.1). By the assumptions, there exists a constant $C > 0$ such that $|K(s) - K(t)| \leq C|s - t|$, $K(s) \leq C$, $|f(s) - f(t)| \leq C|s - t|$ and $|f(s)| \leq C$ for any $s$ and $t$. Define $\mathbb{J} = \{j \in \mathbb{Z} : n(-bL + \tau - d_1)/(d_2 - d_1) \leq j \leq n(bL + \tau - d_1)/(d_2 - d_1)\} = \{l_1, \ldots, l_2\}$. Clearly $\#\mathbb{J} \leq 2nbL/(d_2 - d_1) + 2$, $K_b(u_j - \tau) = 0$ for any $j \notin \mathbb{J}$ and $K_b(u - \tau) = 0$ for $u \leq u_{l_1 - 1}$ or $u \geq u_{l_2 + 1}$. It follows that

$$\left|\frac{1}{n}\sum_{j=1}^{n} K_b(u_j - \tau)f(u_j) - \frac{1}{d_2 - d_1}\int_{d_1}^{d_2} K_b(u - \tau)f(u)\,du\right|$$

$$= \frac{1}{d_2 - d_1}\left|\sum_{j \in \mathbb{J}}\int_{u_{j-1}}^{u_j} K_b(u_j - \tau)f(u_j)\,du - \sum_{j \in \mathbb{J}}\int_{u_{j-1}}^{u_j} K_b(u - \tau)f(u)\,du - \int_{u_{l_2}}^{u_{l_2+1}} K_b(u - \tau)f(u)\,du\right|$$

$$\leq \frac{1}{d_2 - d_1}\sum_{j \in \mathbb{J}}\left\{\int_{u_{j-1}}^{u_j} \frac{C(u_j - u)|f(u_j)|}{b^2}\,du + \int_{u_{j-1}}^{u_j} K_b(u - \tau)C(u_j - u)\,du\right\} + \frac{C^2}{nb}$$

$$\leq \left(\frac{nbL}{d_2 - d_1} + 1\right)\left(\frac{C^2}{b^2} + \frac{C^2}{b}\right)\frac{(d_2 - d_1)}{n^2} + \frac{C^2}{nb} = O\left(\frac{1}{nb}\right).$$

We now show (A.2). Following (A.1),

$$\sup_{t \in [bL, 1-bL]} E\left[\left\{\frac{1}{n}\sum_{j=1}^{n}\epsilon_j K_b(t_j - t) - E(\epsilon_1)\right\}^2\right]$$

$$= \sup_{t \in [bL, 1-bL]} \text{var}\left\{\frac{1}{n}\sum_{j \in \mathbb{J}}\epsilon_j K_b(t_j - t)\right\}$$

$$+ \sup_{t \in [bL,\ldots,1-bL]} \left\{ \frac{1}{n} \sum_{j=1}^{n} K_b(t_j - t) - 1 \right\}^2 \{E(\epsilon_1)\}^2$$

$$\leq \frac{C^2 g}{n^2 b^2}(2nbL + 2)\operatorname{var}(\epsilon_1) + O\left(\frac{1}{n^2 b^2}\right) = O\left(\frac{1}{nb}\right). \qquad \square$$

COROLLARY A.1. *Assume Conditions A2, A5 and A6.*
1. *For any $l = 0, 1, 2, \ldots$, we have that uniformly in $t \in [0,1]$,*

(A.3) $$\frac{1}{n}\sum_{j=1}^{n} K_b(t_j - t)(t_j - t)^l = b^l \left\{ \int_{-t/b}^{(1-t)/b} K(u) u^l \, du + O\left(\frac{1}{nb}\right) \right\},$$

*and thus, uniformly in $i \in [nbL, \ldots, n - nbL]$,*

(A.4) $$\frac{1}{n}\sum_{j=1}^{n} K_b(t_j - t_i)(t_j - t_i)^l = b^l \left\{ \int_{-L}^{L} K(u) u^l \, du + O\left(\frac{1}{nb}\right) \right\}.$$

2. *There exists $C > 0$ such that for all $n = 1, 2, \ldots$, $i \in \{1, \ldots, n\}$ and $j \in \{1, \ldots, n\}$,*

(A.5) $$nb|S_d(i,j)| \leq C.$$

*Moreover, $S_d$ is AVUS. Furthermore, for all $n = 1, 2, \ldots$, $i \in [nbL, \ldots, n-nbL]$ and $j \in \{1, \ldots, n\}$,*

(A.6) $$S_d(i,j) = H(i,j)(1 + c_{n,i}),$$

*where $\sup_{nbL \leq i \leq n-nbL} |c_{n,i}| = O\{1/(nb)\}$.*

3. *Let $\{\epsilon_j\}_{j=1}^{\infty}$ be a sequence of $g$-dependent and identically distributed random variables. Assume $E(\epsilon_1) = 0$ and $E(\epsilon_1^2) < \infty$. Let $Y_i = n^{-1}\sum_{j=1}^{n} \epsilon_j K_b(t_j - t_i)$, $\boldsymbol{\epsilon} = (\epsilon_1, \ldots, \epsilon_n)^T$ and $\mathbf{y} = (Y_1, \ldots, Y_n)^T$. Assume that $B \in \mathbb{R}^{n \times n}$ is AVUS. Then,*

(A.7) $$n^{-1}E(\|B\mathbf{y}\|^2) = o(1),$$

(A.8) $$n^{-1}E(\|BS_d\boldsymbol{\epsilon}\|^2) = o(1).$$

PROOF. Part 1. Following the proof of (A.1),

$$\left| \frac{1}{n}\sum_{j=1}^{n} K_b(t_j - t)\left(\frac{t_j - t}{b}\right)^l - \int_{-t/b}^{(1-t)/b} K(u) u^l \, du \right|$$

$$= \left| \sum_{j \in \mathbb{J}} \int_{t_{j-1}}^{t_j} K_b(t_j - t)\left(\frac{t_j - t}{b}\right)^l du \right.$$

$$\left. - \sum_{j \in \mathbb{J}} \int_{t_{j-1}}^{t_j} K_b(u - t)\left(\frac{u - t}{b}\right)^l du - \int_{t_{l_2}}^{t_{l_2+1}} K_b(u - t)\left(\frac{u - t}{b}\right)^l du \right|$$

$$\leq \sum_{j \in \mathbb{J}} \int_{t_{j-1}}^{t_j} |K_b(t_j - t) - K_b(u - t)| \left|\frac{t_j - t}{b}\right|^l du$$

$$+ \sum_{j \in \mathbb{J}} \int_{t_{j-1}}^{t_j} K_b(u - t) \left|\left(\frac{t_j - t}{b}\right)^l - \left(\frac{u - t}{b}\right)^l\right| du + \frac{1}{nb}.$$

Note that for $j \in \mathbb{J}$, we have $|(t_j - t)/b|^l \leq L^l$ and that for $j \in \mathbb{J}$ and $t_{j-1} \leq u \leq t_j$, we have $|\{(t_j - t)/b\}^l - \{(u - t)/b\}^l| \leq C|(t_j - u)/b|$. Applying the same argument for (A.1) completes the proof for (A.3) and, in turn, (A.4).

We then show Part 2. Define $S_{n,l}(t) = n^{-1} \sum_{j=1}^n K_b(t_j - t)(t_j - t)^l$, $l = 0, 1, 2$, and $S_n(t) = \mathbf{X}(t)^T \mathbf{W}(t) \mathbf{X}(t)$. Then,

$$S_n(t) = n \begin{bmatrix} S_{n,0}(t) & S_{n,1}(t) \\ S_{n,1}(t) & S_{n,2}(t) \end{bmatrix},$$

$$\{S_n(t)\}^{-1} = \frac{1}{n[S_{n,0}(t)S_{n,2}(t) - \{S_{n,1}(t)\}^2]} \begin{bmatrix} S_{n,2}(t) & -S_{n,1}(t) \\ -S_{n,1}(t) & S_{n,0}(t) \end{bmatrix}.$$

According to (A.3), uniformly in $t \in [0, 1]$,

$$\{S_n(t)\}^{-1} = \frac{1}{nb^2[f(t) + O\{1/(nb)\}]}$$
$$\times \begin{bmatrix} b^2[a_2(t) + O\{1/(nb)\}] & -b[a_1(t) + O\{1/(nb)\}] \\ -b[a_1(t) + O\{1/(nb)\}] & a_0(t) + O\{1/(nb)\} \end{bmatrix},$$

where $a_l(t) = \int_{-t/b}^{(1-t)/b} K(u) u^l \, du$, $l = 0, 1, 2$, are all uniformly bounded in $t$ and $f(t) = a_0(t)a_2(t) - \{a_1(t)\}^2$ is minimized at $t = 0$ with $f(0) = 0.25 \, \text{var}(|U|) > 0$ for a random variable $U$ with density $K(u)$. It is seen from the definition of $S_d$ in (3.3) that

$$S_d(i, j) = \frac{K_b(t_j - t_i)/n}{f(t_i) + O\{1/(nb)\}} \left[\left\{a_2(t_i) + O\left(\frac{1}{nb}\right)\right\} - \frac{t_j - t_i}{b} \left\{a_1(t_i) + O\left(\frac{1}{nb}\right)\right\}\right].$$

Note that when $j \notin [i - nbL, \ldots, i + nbL]$, $K_b(t_j - t_i) = 0$ implies $S_d(i, j) = 0$. Also, note that when $j \in [i - nbL, \ldots, i + nbL]$, $|t_j - t_i|/b \leq L$. Thus, there exists $C > 0$ uniformly in $i \in \{1, \ldots, n\}$ and $j \in [i - nbL, \ldots, i + nbL]$ such that $|S_d(i, j)| \leq Cn^{-1} K_b(t_j - t_i) \leq C \sup_t K(t)/(nb)$. Thus, for some $C > 0$, uniformly in $n = 1, 2, \ldots$, $i = 1, \ldots, n$ and $j = 1, \ldots, n$,

$$|S_d(i, j)| \begin{cases} = 0, & \text{if } |j - i| > nbL, \\ \leq C/(nb), & \text{if } |j - i| \leq nbL, \end{cases}$$

which implies (A.5) and also implies that $S_d$ is AVUS. Similar arguments for (A.5) combined with (A.4) and (3.3) yield (A.6).

We now show Part 3. From (A.2), for any $\epsilon > 0$, there exists $N$ such that $n > N$ implies that $E(Y_i^2) < \epsilon$ for all $i \in [nbL, \ldots, n - nbL]$. It can also be shown that there exists $C > 0$ such that $E(Y_i^2) \leq C$ for all $i \in [1, \ldots, nbL] \cup [n - nbL, \ldots, n]$. Since $B$ is AVUS, there exists $C_1 > 0$ such that $\sup_{n \geq 1} \sup_{1 \leq j \leq n} \sum_{i=1}^n |B(i,j)| \leq C_1$ and $\sup_{n \geq 1} \sup_{1 \leq i \leq n} \sum_{j=1}^n |B(i,j)| \leq C_1$. The proof of (A.7) is obtained as follows:

$$n^{-1} E(\|B\mathbf{y}\|^2) = \frac{1}{n} \sum_{i=1}^n \sum_{j=1}^n \sum_{k=1}^n B(i,j) B(i,k) E(Y_j Y_k)$$

$$\leq \frac{1}{n} \sum_{i=1}^n \sum_{j=1}^n \sum_{k=1}^n |B(i,j) B(i,k)| \{E(Y_j^2) E(Y_k^2)\}^{1/2}$$

$$\leq \epsilon^{1/2} C^{1/2} \frac{1}{n} \sum_{i=1}^n \left\{ \sum_{j=1}^n |B(i,j)| \right\}^2$$

$$+ \frac{C}{n} \sum_{i=1}^n \sum_{j,k \in [1,\ldots,nbL] \cup [n-nbL,\ldots,n]} |B(i,j) B(i,k)|$$

$$\leq \epsilon^{1/2} C^{1/2} C_1^2 + 2C_1^2 CbL + 2C_1^2 C/n.$$

Applying (A.5), (A.6) and similar arguments for (A.7) completes the proof of (A.8). □

LEMMA A.4. *Let $\{X_n\}_{n=1}^\infty$ be a sequence of random variables such that every subsequence of $X_n$ has a further subsequence converging in distribution to a same random variable $X$. Then $X_n \xrightarrow{\mathcal{L}} X$.*

PROOF. Let $\phi_Y$ denote the characteristic function of a random variable $Y$. Since any subsequence $\{n_l\}_{l=1}^\infty$ of $\{1, 2, \ldots\}$ has a further subsequence, $\{n_{l_j}\}_{j=1}^\infty$, such that $X_{n_{l_j}} \xrightarrow{\mathcal{L}} X$, the Lévy–Cramér continuity theorem [Shao (2003), page 56] implies that $\phi_{X_{n_{l_j}}}(t) \to \phi_X(t)$ as $j \to \infty$ for any $t \in \mathbb{R}$. This, in turn indicates that $\phi_{X_n}(t) \to \phi_X(t)$ as $n \to \infty$ for any $t \in \mathbb{R}$. We then conclude that $X_n \xrightarrow{\mathcal{L}} X$ by repeated application of the Lévy–Cramér continuity theorem. □

LEMMA A.5. *Let $\{\epsilon_i\}_{i=1}^n$ be a stationary g-dependent sequence with $E(\epsilon_1) = 0$ and $E(\epsilon_1^4) < \infty$. Set $x_{n,i} = \tau_{n,i} \epsilon_i$, $i = 1, \ldots, n$, where $\{\tau_{n,i}\}$ is independent of $\{\epsilon_i\}$. Define $\sigma_n^2(\{\tau_{n,i}\}) = E[(\sum_{i=1}^n x_{n,i})^2 | \{\tau_{n,i}\}]$. If* $\sup_{n\geq 1}\sup_{1\leq i\leq n}|\tau_{n,i}|\leq C$ and $\sigma_n^2(\{\tau_{n,i}\})=n\sigma_a^2\{1+o_P(1)\}$ *for some constants* $C>0$ *and* $\sigma_a^2\in(0,\infty)$, *then* $\sigma_n^{-1}\sum_{i=1}^n x_{n,i}\xrightarrow{\mathcal{L}} N(0,1)$.

PROOF. The proof follows directly from Theorem 1.7 of Bosq (1998), page 36. It can be achieved by applying the blocking arguments and Lyapunov's central limit theorem. □

LEMMA A.6. *Assume model* (3.2) *and Conditions* A1–A7. *Then:*

1. *all three matrices* $\widetilde{V}_R = V(\mathbf{I}-S_d)$, $\widetilde{V}_L = (\mathbf{I}-S_d)^T V$ *and* $\widetilde{V} = (\mathbf{I}-S_d)^T V(\mathbf{I}-S_d)$ *are* AVUS;
2. $\mathrm{var}(n^{-1}\boldsymbol{\xi}_{j_1,\ell_1}^T\widetilde{V}\boldsymbol{\xi}_{j_2,\ell_2})\to 0$ *for all* $j_1,j_2=1,\ldots,r$ *and all* $\ell_1,\ell_2=1,\ldots,m$;
3. $n^{-1}\widetilde{\mathbf{S}}^T R^{-1}\widetilde{\mathbf{S}} - E(n^{-1}\widetilde{\mathbf{S}}^T R^{-1}\widetilde{\mathbf{S}})\xrightarrow{\mathcal{P}}\mathbf{0}$;
4. *all entries of* $E(n^{-1}\widetilde{\mathbf{S}}^T R^{-1}\widetilde{\mathbf{S}})$ *are bounded;*
5. *all convergent subsequences of* $E(n^{-1}\widetilde{\mathbf{S}}^T R^{-1}\widetilde{\mathbf{S}})$ *are positive definite.*

PROOF. The proof of Part 1 can be obtained from applying Lemma A.1, Lemma A.2 and part 2 of Corollary A.1.

We next show part 2. For $\ell_1,\ell_2=1,\ldots,m$,

$$n^{-1}\boldsymbol{\xi}_{j_1,\ell_1}^T\widetilde{V}\boldsymbol{\xi}_{j_2,\ell_2} = n^{-1}\sum_{k_1=1}^n\sum_{k_2=1}^n \mathbf{S}_{j_1}(k_1,\ell_1)\widetilde{V}(k_1,k_2)\mathbf{S}_{j_2}(k_2,\ell_2)$$

$$= n^{-1}\sum_{k_1=\ell_1}^n\sum_{k_2=\ell_2}^n s_{j_1}(t_{k_1}-t_{\ell_1})s_{j_2}(t_{k_2}-t_{\ell_2})\widetilde{V}(k_1,k_2)$$

$$\equiv I_{1,1}+I_{1,2}+I_{1,3}+I_{1,4},$$

where

$$I_{1,1}=n^{-1}p_{j_1}p_{j_2}\sum_{k_1=\ell_1}^n\sum_{k_2=\ell_2}^n \widetilde{V}(k_1,k_2),$$

$$I_{1,2}=n^{-1}p_{j_2}\sum_{k_1=\ell_1}^n\{s_{j_1}(t_{k_1}-t_{\ell_1})-p_{j_1}\}\sum_{k_2=\ell_2}^n \widetilde{V}(k_1,k_2),$$

$$I_{1,3}=n^{-1}p_{j_1}\sum_{k_2=\ell_2}^n\{s_{j_2}(t_{k_2}-t_{\ell_2})-p_{j_2}\}\sum_{k_1=\ell_1}^n \widetilde{V}(k_1,k_2),$$

$$I_{1,4}=n^{-1}\sum_{k_1=\ell_1}^n\sum_{k_2=\ell_2}^n\{s_{j_1}(t_{k_1}-t_{\ell_1})-p_{j_1}\}\{s_{j_2}(t_{k_2}-t_{\ell_2})-p_{j_2}\}\widetilde{V}(k_1,k_2).$$

It is easily seen that

(A.9) $\mathrm{var}(I_{1,1})=0,\quad \mathrm{var}(I_{1,2})=O(n^{-1})\quad\text{and}\quad \mathrm{var}(I_{1,3})=O(n^{-1}).$

For $I_{1,4}$,

$$E(I_{1,4}^2)$$
$$= \frac{1}{n^2} \sum_{k_1=\ell_1}^{n} \sum_{k_2=\ell_2}^{n} \sum_{k_3=\ell_1}^{n} \sum_{k_4=\ell_2}^{n} E[\{s_{j_1}(t_{k_1} - t_{\ell_1}) - p_{j_1}\}\{s_{j_2}(t_{k_2} - t_{\ell_2}) - p_{j_2}\}$$
$$\times \{s_{j_1}(t_{k_3} - t_{\ell_1}) - p_{j_1}\}\{s_{j_2}(t_{k_4} - t_{\ell_2}) - p_{j_2}\}]$$
$$\times \widetilde{V}(k_1, k_2)\widetilde{V}(k_3, k_4),$$

in which two situations will be discussed. In the situation where $j_1 = j_2 = j$, the additive term above is nonzero only in the following four cases:

I:    $k_1 - \ell_1 = k_2 - \ell_2 = k_3 - \ell_1 = k_4 - \ell_2$;
II:    $\{k_1 - \ell_1 = k_2 - \ell_2\} \neq \{k_3 - \ell_1 = k_4 - \ell_2\}$;
III:    $\{k_1 - \ell_1 = k_3 - \ell_1\} \neq \{k_2 - \ell_2 = k_4 - \ell_2\}$;
IV:    $\{k_1 - \ell_1 = k_4 - \ell_2\} \neq \{k_2 - \ell_2 = k_3 - \ell_1\}$.

Thus, $E(I_{1,4}^2) = E_{\text{I}} + E_{\text{II}} + E_{\text{III}} + E_{\text{IV}}$, where

$$E_{\text{I}} \leq n^{-2}[p_j(1 - p_j)\{p_j^3 + (1 - p_j)^3\}] \sum_{k_1=1}^{n} \{\widetilde{V}(k_1, k_1 + \ell_2 - \ell_1)\}^2 = O(n^{-2}),$$

$$E_{\text{II}} \leq n^{-2}[p_j(1 - p_j)]^2 \sum_{k_1=1}^{n} \sum_{k_3=1}^{n} |\widetilde{V}(k_1, k_1 + \ell_2 - \ell_1)| \cdot |\widetilde{V}(k_3, k_3 + \ell_2 - \ell_1)|$$
$$= O(n^{-2}),$$

$$E_{\text{III}} \leq n^{-2}[p_j(1 - p_j)]^2 \sum_{k_1=1}^{n} \sum_{k_2=1}^{n} \{\widetilde{V}(k_1, k_2)\}^2 = O(n^{-2}),$$

$$E_{\text{IV}} \leq n^{-2}[p_j(1 - p_j)]^2 \sum_{k_1=1}^{n} \sum_{k_2=1}^{n} |\widetilde{V}(k_1, k_2)| \cdot |\widetilde{V}(k_2 - \ell_2 + \ell_1, k_1 + \ell_2 - \ell_1)|$$
$$= O(n^{-2}).$$

Hence, $E(I_{1,4}^2) = O(n^{-2})$ when $j_1 = j_2$. In the situation where $j_1 \neq j_2$, since the $s_{j_1}(\cdot)$ are independent at different time points and, similarly, the $s_{j_2}(\cdot)$ are independent at different time points, $E[\{s_{j_1}(t_{k_1} - t_{\ell_1}) - p_{j_1}\}\{s_{j_1}(t_{k_3} - t_{\ell_1}) - p_{j_1}\}\{s_{j_2}(t_{k_2} - t_{\ell_2}) - p_{j_2}\}\{s_{j_2}(t_{k_4} - t_{\ell_2}) - p_{j_2}\}]$ is nonzero only if $k_1 = k_3$ and $k_2 = k_4$. In this case,

$$E(I_{1,4}^2)$$
$$= \frac{1}{n^2} \sum_{k_1=\ell_1}^{n} \sum_{k_3=\ell_1}^{n} \sum_{k_2=\ell_2}^{n} \sum_{k_4=\ell_2}^{n} E[\{s_{j_1}(t_{k_1} - t_{\ell_1}) - p_{j_1}\}\{s_{j_1}(t_{k_3} - t_{\ell_1}) - p_{j_1}\}$$
$$\times \{s_{j_2}(t_{k_2} - t_{\ell_2}) - p_{j_2}\}\{s_{j_2}(t_{k_4} - t_{\ell_2}) - p_{j_2}\}]$$

$$\times \widetilde{V}(k_1, k_2)\widetilde{V}(k_3, k_4)$$

$$= \frac{1}{n^2} \sum_{k_1=\ell_1}^{n} \sum_{k_2=\ell_2}^{n} E[\{s_{j_1}(t_{k_1} - t_{\ell_1}) - p_{j_1}\}^2 \{s_{j_2}(t_{k_2} - t_{\ell_2}) - p_{j_2}\}^2]$$

$$\times \{\widetilde{V}(k_1, k_2)\}^2$$

$$= n^{-2} C \sum_{k_1=\ell_1}^{n} \sum_{k_2=\ell_2}^{n} \{\widetilde{V}(k_1, k_2)\}^2 = O(n^{-2}).$$

Thus, in both situations, $\mathrm{var}(I_{1,4}) \to 0$. This, combined with (A.9), yields Part 2.

We then show Part 3. Recall that $\widetilde{\mathbf{S}} = (\mathbf{I} - S_d)\mathbf{S} = [\widetilde{\mathbf{S}}_1, \ldots, \widetilde{\mathbf{S}}_r]$, where $\widetilde{\mathbf{S}}_j = (\mathbf{I} - S_d)\mathbf{S}_j$. Then,

$$n^{-1}\widetilde{\mathbf{S}}^T R^{-1} \widetilde{\mathbf{S}} = n^{-1} \begin{bmatrix} \widetilde{\mathbf{S}}_1^T R^{-1} \widetilde{\mathbf{S}}_1 & \cdots & \widetilde{\mathbf{S}}_1^T R^{-1} \widetilde{\mathbf{S}}_r \\ \vdots & \ddots & \vdots \\ \widetilde{\mathbf{S}}_r^T R^{-1} \widetilde{\mathbf{S}}_1 & \cdots & \widetilde{\mathbf{S}}_r^T R^{-1} \widetilde{\mathbf{S}}_r \end{bmatrix}.$$

It suffices to consider the block matrix $n^{-1}\widetilde{\mathbf{S}}_{j_1}^T R^{-1} \widetilde{\mathbf{S}}_{j_2}$, whose $(\ell_1, \ell_2)$th entry is $C_{j_1,\ell_1;j_2,\ell_2} = n^{-1}\boldsymbol{\xi}_{j_1,\ell_1}^T \widetilde{V} \boldsymbol{\xi}_{j_2,\ell_2}$. By Part 2, $\mathrm{var}(C_{j_1,\ell_1;j_2,\ell_2}) \to 0$, which, in turn, gives $C_{j_1,\ell_1;j_2,\ell_2} - E(C_{j_1,\ell_1;j_2,\ell_2}) \xrightarrow{\mathcal{P}} 0$ and the conclusion of Part 3.

We now show Part 4, which can easily be derived from

$$|n^{-1}\boldsymbol{\xi}_{j_1,\ell_1}^T \widetilde{V} \boldsymbol{\xi}_{j_2,\ell_2}| \leq n^{-1} \sum_{k_1=1}^{n} \sum_{k_2=1}^{n} |\widetilde{V}(k_1, k_2)| \leq C$$

since the entries of $\mathbf{S}_j$ are either 0 or 1.

Last, we show Part 5. For any $\{n_k\}_{k=1}^{\infty}$ such that $E(n_k^{-1}\widetilde{\mathbf{S}}^T R^{-1} \widetilde{\mathbf{S}})$ converges to some limit $\mathbf{M}$, by Part 3, $n_k^{-1}\widetilde{\mathbf{S}}^T R^{-1} \widetilde{\mathbf{S}} \xrightarrow{\mathcal{P}} \mathbf{M}$. Obviously, $\mathbf{M}$ is semi-positive definite. It remains to show that $\mathbf{M}$ is nonsingular. We now prove this by contradiction. Assume that there exists some $\mathbf{c} = (\mathbf{c}_1^T, \ldots, \mathbf{c}_r^T)^T \in \mathbb{R}^{rm}$, where $\mathbf{c}_j = (c_{j,1}, \ldots, c_{j,m})^T$, $j = 1, \ldots, r$, such that $\mathbf{c} \neq \mathbf{0}$ and $\mathbf{c}^T \mathbf{M} \mathbf{c} = 0$. Then, $n_k^{-1}\mathbf{c}^T \widetilde{\mathbf{S}}^T R^{-1} \widetilde{\mathbf{S}} \mathbf{c} \xrightarrow{\mathcal{P}} 0$. By the Schur decomposition, there exist $Q$ and $v_1, \ldots, v_{n_k}$ such that $V = Q^T \mathrm{diag}(v_1, \ldots, v_{n_k}) Q$, where $Q^T Q = \mathbf{I}_{n_k}$ and $v_{n_k} \leq \cdots \leq v_1$. Furthermore, from Condition A3, $V$ is positive definite with eigenvalues bounded away from 0 and $\infty$. Thus, there exist constants $a$ and $b$ such that $0 < a < b < \infty$ and $0 < a \leq v_{n_k} \leq \cdots \leq v_1 \leq b < \infty$. Noting that

$$n_k^{-1}\mathbf{c}^T \widetilde{\mathbf{S}}^T R^{-1} \widetilde{\mathbf{S}} \mathbf{c} = \widetilde{\mathbf{c}}^T \mathrm{diag}(v_1, \ldots, v_{n_k}) \widetilde{\mathbf{c}} \geq a \|\widetilde{\mathbf{c}}\|^2,$$

where $\widetilde{\mathbf{c}} \equiv Q\widetilde{\mathbf{S}}\mathbf{c}/\sqrt{n_k}$, we conclude that as $k \to \infty$,

(A.10) $$\|(\mathbf{I} - S_d)\mathbf{S}\mathbf{c}\|^2 / n_k \xrightarrow{\mathcal{P}} 0.$$

Consider

$$\|(\mathbf{I} - S_d)\mathbf{Sc}\|^2/n = \|\{\mathbf{Sc} - E(\mathbf{S})\mathbf{c}\} - \{S_d\mathbf{Sc} - E(\mathbf{S})\mathbf{c}\}\|^2/n$$
$$\equiv J_1 - 2J_2 + J_3, \quad \text{(A.11)}$$

where $J_1 = \|\mathbf{Sc} - E(\mathbf{S})\mathbf{c}\|^2/n$, $J_2 = \{\mathbf{Sc} - E(\mathbf{S})\mathbf{c}\}^T\{S_d\mathbf{Sc} - E(\mathbf{S})\mathbf{c}\}/n$ and $J_3 = \|S_d\mathbf{Sc} - E(\mathbf{S})\mathbf{c}\|^2/n$. By the law of large numbers and block arguments, we can show that

$$\text{(A.12)} \quad J_1 \xrightarrow{\mathcal{P}} \text{var}\left(\sum_{j=1}^r \mathbf{c}_j^T \mathbf{s}_{j,m}\right) \geq \left(1 - \sum_{j=1}^r p_j\right)\left(\sum_{j=1}^r p_j\|\mathbf{c}_j\|^2\right) > 0,$$

where $\mathbf{s}_{j,i} = \mathbf{S}_j^T \mathbf{e}_i$, $i = 1, \ldots, m$. For $J_3$, note that

$$E(J_3) \leq o(1) + \frac{\|\mathbf{c}\|^2}{n} \sum_{i=1}^n E\left[\sum_{j=1}^r \sum_{k=1}^m \left\{\sum_{\ell=1}^n S_d(i,\ell) s_j(t_\ell - t_k) - p_j\right\}^2\right]$$

$$\leq o(1) + \frac{\|\mathbf{c}\|^2}{n} \sum_{j=1}^r \sum_{k=1}^m (n - 2nbL)$$

$$\times \sup_{i \in [nbL, \ldots, n-nbL]} E\left[\left\{\sum_{\ell=1}^n S_d(i,\ell) s_j(t_\ell - t_k) - p_j\right\}^2\right]$$

$$+ \frac{\|\mathbf{c}\|^2}{n} \sum_{j=1}^r \sum_{k=1}^m 2nbL$$

$$\times \sup_{i \in [1, \ldots, nbL] \cup [n-nbL, \ldots, n]} E\left[\left\{\sum_{\ell=1}^n S_d(i,\ell) s_j(t_\ell - t_k) - p_j\right\}^2\right]$$

$$= o(1) + \frac{\|\mathbf{c}\|^2}{n} \sum_{j=1}^r \sum_{k=1}^m (n - 2nbL)$$

$$\times \sup_{i \in [nbL, \ldots, n-nbL]} E\left[\left\{\{1 + o(1)\}\frac{1}{n}\sum_{\ell=1}^n K_b(t_\ell - t_i) s_j(t_\ell - t_k) - p_j\right\}^2\right]$$

$$+ \frac{\|\mathbf{c}\|^2}{n} \sum_{j=1}^r \sum_{k=1}^m 2nbL$$

$$\times \sup_{i \in [1, \ldots, nbL] \cup [n-nbL, \ldots, n]} E\left[\left\{\sum_{\ell=1}^n S_d(i,\ell) s_j(t_\ell - t_k) - p_j\right\}^2\right]$$

$$= o(1) + \|\mathbf{c}\|^2 \left\{\sum_{j=1}^r \sum_{k=1}^m o(1)\right\} + \|\mathbf{c}\|^2 \left\{\sum_{j=1}^r \sum_{k=1}^m O(b)\right\}.$$

In the last equality, the second term follows from (A.2), whereas the third term uses $|s_j(t_\ell - t_k)| \leq 1$ and the fact that $S_d$ is AVUS. Thus, $E(J_3) \to 0$ and $J_3 \geq 0$ imply that

$$J_3 = o_P(1). \tag{A.13}$$

By the Cauchy–Schwarz inequality, $|J_2| \leq 2(J_1 J_3)^{1/2}$. Thus, $J_2 = o_P(1)$. This, together with (A.11), (A.12) and (A.13), shows that $\|(\mathbf{I} - S_d)\mathbf{Sc}\|^2/n \xrightarrow{\mathcal{P}} \mathrm{var}(\sum_{j=1}^r \mathbf{c}_j^T \mathbf{s}_{j,m})$, which contradicts (A.10). □

LEMMA A. 7. *Assume Condition* A. *Suppose that* $n^{-1}\widetilde{\mathbf{S}}^T R^{-1} \widetilde{\mathbf{S}} \xrightarrow{\mathcal{P}} \mathbf{M}$, *where* $\mathbf{M} \in \mathbb{R}^{rm \times rm}$ *is positive definite. Then* $n^{-1}\widetilde{\mathbf{S}}^T \widehat{R}^{-1} \widetilde{\mathbf{S}} \xrightarrow{\mathcal{P}} \mathbf{M}$ *and* $n^{1/2}(\widehat{\mathbf{h}} - \mathbf{h}) \xrightarrow{\mathcal{L}} N(\mathbf{0}, \sigma^2 \mathbf{M}^{-1})$.

PROOF. From (3.5), $\widehat{\mathbf{h}} - \mathbf{h} = (n^{-1}\widetilde{\mathbf{S}}^T \widehat{R}^{-1} \widetilde{\mathbf{S}})^{-1} n^{-1/2}(J_1^* + J_2^*)$, where $J_1^* = n^{-1/2}\widetilde{\mathbf{S}}^T \widehat{R}^{-1} \widetilde{\mathbf{d}}$ and $J_2^* = n^{-1/2}\widetilde{\mathbf{S}}^T \widehat{R}^{-1} \widetilde{\boldsymbol{\epsilon}}$. Let $J_1 = n^{-1/2}\widetilde{\mathbf{S}}^T R^{-1} \widetilde{\mathbf{d}}$ and $J_2 = n^{-1/2}\widetilde{\mathbf{S}}^T R^{-1} \widetilde{\boldsymbol{\epsilon}}$. The proof proceeds in three steps to show that $J_1 = o_P(1)$, $J_2 \xrightarrow{\mathcal{L}} N(\mathbf{0}, \sigma^2 \mathbf{M})$ and

$$\begin{aligned} n^{-1}\widetilde{\mathbf{S}}^T(\widehat{R}^{-1} - R^{-1})\widetilde{\mathbf{S}} &\xrightarrow{\mathcal{P}} \mathbf{0}, \\ J_1^* - J_1 &= o_P(1), \\ J_2^* - J_2 &= o_P(1). \end{aligned} \tag{A.14}$$

First, we will show $J_1 = o_P(1)$. From (A.6), the $i$th entry of $(\mathbf{I} - S_d)\mathbf{d}$ is

$$d(t_i) - \{1 + o(1)\}\frac{1}{n}\sum_{j=1}^n K_b(t_j - t_i)\,d(t_j)$$

$$= d(t_i) - \{1 + o(1)\} \int_{-t_i/b}^{(1-t_i)/b} K(u)\,d(t_i + ub)\,du + O\left(\frac{1}{nb}\right)$$

$$= \{1 + o(1)\} \int_{-L}^{L} \{ubK(u)d'(t_i) + u^2 b^2 K(u) d''(\xi)/2\}\,du + o(1) + O\left(\frac{1}{nb}\right)$$

$$= o(1),$$

uniformly in $i \in [nbL, \ldots, n - nbL]$. When $i \in [1, \ldots, nbL] \cup [n - nbL, \ldots, n]$,

$$\left|d(t_i) - \{1 + o(1)\}\frac{1}{n}\sum_{j=1}^n K_b(t_j - t_i)d(t_j)\right|$$

$$\leq \sup_{t \in [0,1]} |d(t)| \left\{1 + \sup_{1 \leq \ell \leq n}\left|\frac{1}{n}\sum_{j=1}^n K_b(t_j - t_\ell)\right|\right\} \leq C,$$

for some $C > 0$. Thus,

(A.15)
$$\sup_{i \in [nbL, \ldots, n-nbL]} |e_i^T \widetilde{\mathbf{d}}| = o(1),$$
$$\sup_{n \geq 1} \sup_{i \in [1, \ldots, nbL] \cup [n-nbL, \ldots, n]} |e_i^T \widetilde{\mathbf{d}}| \leq C.$$

Consider the $j$th block vector of $J_1$: $J_{1,j} = n^{-1/2} \widetilde{\mathbf{S}}_j^T R^{-1} \widetilde{\mathbf{d}}$. Its $i$th entry is $e_i^T J_{1,j} = n^{-1/2} \widetilde{\mathbf{d}}^T V(\mathbf{I} - S_d)(\boldsymbol{\xi}_{j,i} - p_j \mathbf{1})$. Then,

$$E(e_i^T J_{1,j}) = n^{-1/2} \widetilde{\mathbf{d}}^T V(\mathbf{I} - S_d)\{p_j (\mathbf{0}_{i-1}^T, \mathbf{1}_{n-i+1}^T)^T - p_j \mathbf{1}\}$$
$$= -p_j n^{-1/2} \widetilde{\mathbf{d}}^T V(\mathbf{I} - S_d)(\mathbf{1}_{i-1}^T, \mathbf{0}_{n-i+1}^T)^T$$
$$= -p_j n^{-1/2} \sum_{k=1}^n \widetilde{d}(t_k) \sum_{\ell=1}^{i-1} \widetilde{V}_R(k, \ell),$$

and thus $|E(e_i^T J_{1,j})| \leq p_j n^{-1/2} \sum_{\ell=1}^m \sum_{k=1}^n |\widetilde{V}_R(k, \ell)| \{\sup_{1 \leq l \leq n} |\widetilde{d}(t_l)|\} = o(1)$, by (A.15) and the fact that $\widetilde{V}_R$ is AVUS. Moreover,

$$\text{var}(e_i^T J_{1,j}) = n^{-1} \widetilde{\mathbf{d}}^T V(\mathbf{I} - S_d) \begin{bmatrix} \mathbf{0} & \mathbf{0} \\ \mathbf{0} & p_j(1 - p_j)\mathbf{I} \end{bmatrix} (\mathbf{I} - S_d)^T V \widetilde{\mathbf{d}}$$
$$\leq n^{-1} p_j (1 - p_j) \|\widetilde{V}_L \widetilde{\mathbf{d}}\|^2,$$

which implies that $\text{var}(J_{1,j}) = \mathbf{1}_m \mathbf{1}_m^T o(1)$, using similar derivations for (A.7). Thus, $J_{1,j} \xrightarrow{\mathcal{P}} \mathbf{0}$ and hence $J_1 \xrightarrow{\mathcal{P}} \mathbf{0}$.

Second, we will show that $J_2 \xrightarrow{\mathcal{L}} N(\mathbf{0}, \sigma^2 \mathbf{M})$. Since $\widetilde{\boldsymbol{\epsilon}} = \boldsymbol{\epsilon} - S_d \boldsymbol{\epsilon}$, $J_2 = n^{-1/2} \times \widetilde{\mathbf{S}}^T R^{-1} \boldsymbol{\epsilon} - n^{-1/2} \widetilde{\mathbf{S}}^T R^{-1} S_d \boldsymbol{\epsilon} \equiv J_{21} - J_{22}$. Consider the $j$th block vector of $J_{22}$: $J_{22,j} = n^{-1/2} \widetilde{\mathbf{S}}_j^T R^{-1} S_d \boldsymbol{\epsilon}$. Its $i$th entry is $e_i^T J_{22,j} = n^{-1/2} (S_d \boldsymbol{\epsilon})^T V(\mathbf{I} - S_d) \boldsymbol{\xi}_{j,i}$, thus

$$E(e_i^T J_{22,j}) = E[E\{n^{-1/2}(S_d \boldsymbol{\epsilon})^T V(\mathbf{I} - S_d) \boldsymbol{\xi}_{j,i} | \boldsymbol{\xi}_{j,i}\}] = 0$$

and

$$\text{var}(e_i^T J_{22,j}) = \text{var}\{E(e_i^T J_{22,j} | \boldsymbol{\epsilon})\} + E\{\text{var}(e_i^T J_{22,j} | \boldsymbol{\epsilon})\}$$
$$= n^{-1} p_j^2 \text{var}\{(S_d \boldsymbol{\epsilon})^T V(\mathbf{I} - S_d)(\mathbf{1}_{i-1}^T, \mathbf{0}_{n-i+1}^T)^T\}$$
$$+ n^{-1} E\left\{(S_d \boldsymbol{\epsilon})^T V(\mathbf{I} - S_d) \begin{bmatrix} \mathbf{0} & \mathbf{0} \\ \mathbf{0} & p_j(1 - p_j)\mathbf{I} \end{bmatrix} (\mathbf{I} - S_d)^T V(S_d \boldsymbol{\epsilon})\right\}$$
$$\leq n^{-1} p_j^2 \sigma^2 \|R^{1/2} S_d^T V(\mathbf{I} - S_d)(\mathbf{1}_{i-1}^T, \mathbf{0}_{n-i+1}^T)^T\|^2$$
$$+ n^{-1} p_j (1 - p_j) E(\|\widetilde{V}_L S_d \boldsymbol{\epsilon}\|^2)$$
$$= o(1) + o(1) = o(1).$$

In the last equality, the first $o(1)$ is from Lemma A.2 and similar arguments for (A.8). The second $o(1)$ is from (A.8) and Part 1 of Lemma A.6. Thus, $J_{22,j} = o_P(1)$ and $J_{22} = o_P(1)$. For $J_{21}$, by the Cramér–Wold device, it suffices to show that for any $\mathbf{w} = (\mathbf{w}_1^T, \ldots, \mathbf{w}_r^T)^T \in \mathbb{R}^{rm}$, $\mathbf{w}^T J_{21} \xrightarrow{\mathcal{L}} N(0, \sigma^2 \mathbf{w}^T \mathbf{M} \mathbf{w})$, where $\mathbf{w}_j = (w_{j,1}, \ldots, w_{j,m})^T$, $j = 1, \ldots, r$. Note that $\mathbf{w}^T J_{21} = n^{-1/2} \sum_{i=1}^n \tau_{n,i} \epsilon_i$, where $\tau_{n,i} = \sum_{j=1}^r \sum_{k=1}^n \sum_{\ell=1}^m w_{j,\ell} s_j(t_k - t_\ell) \widetilde{V}_L(k,i)$. Thus,

$$|\tau_{n,i}| \leq rm \left( \max_{1 \leq j \leq r} \max_{1 \leq \ell \leq m} |w_{j,\ell}| \right) \sum_{k=1}^n |\widetilde{V}_L(k,i)| \leq C,$$

where the last inequality is from Part 1 of Lemma A.6. Also, $\sigma_n^2(\{\tau_{n,i}\}) = n \operatorname{var}(\mathbf{w}^T J_{21} | \{\tau_{n,i}\}) = \sigma^2 \mathbf{w}^T \widetilde{\mathbf{S}}^T R^{-1} \widetilde{\mathbf{S}} \mathbf{w} = n\sigma^2 \mathbf{w}^T \mathbf{M} \mathbf{w} \{1 + o_P(1)\} = n\sigma_a^2 \{1 + o_P(1)\}$, where $\sigma_a^2 = \sigma^2 \mathbf{w}^T \mathbf{M} \mathbf{w}$. By Lemma A.5, the result follows.

Third, to verify (A.14), it is sufficient to show that

(A.16) $$n^{-1} \widetilde{\mathbf{S}}^T (\widehat{R}^{-1} - R^{-1}) \widetilde{\mathbf{S}} = o_P(1),$$

(A.17) $$n^{-1/2} \widetilde{\mathbf{S}}^T (\widehat{R}^{-1} - R^{-1}) \widetilde{\mathbf{d}} = o_P(1),$$

(A.18) $$n^{-1/2} \widetilde{\mathbf{S}}^T (\widehat{R}^{-1} - R^{-1}) \widetilde{\boldsymbol{\epsilon}} = o_P(1).$$

Note that Condition A8 implies that $\widehat{R}^{-1} - R^{-1}$ is AVUS. Similar arguments for Lemma A.6, $J_1$ and $J_2$ complete the proofs for (A.16), (A.17) and (A.18), respectively. □

COROLLARY A.2. *Assume Condition A. Then,*
1. $\widehat{\mathbf{h}} \xrightarrow{\mathcal{P}} \mathbf{h}$;
2. $(A\widehat{\mathbf{h}} - A\mathbf{h})^T \{A(\widetilde{\mathbf{S}}^T \widehat{R}^{-1} \widetilde{\mathbf{S}})^{-1} A^T\}^{-1} (A\widehat{\mathbf{h}} - A\mathbf{h}) \xrightarrow{\mathcal{L}} \sigma^2 \chi_k^2$.

PROOF. By Lemma A.6, for any subsequence $\{n_l\}_{l=1}^\infty$, there exists a further subsequence, $\{n_{l_j}\}_{j=1}^\infty$, such that $n_{l_j}^{-1} \widetilde{\mathbf{S}}^T R^{-1} \widetilde{\mathbf{S}} \xrightarrow{\mathcal{P}} \mathbf{M}_l$ for some positive definite matrix $\mathbf{M}_l$. For this $\{n_{l_j}\}_{j=1}^\infty$, an appeal to Lemma A.7 gives $n_{l_j}^{-1} \widetilde{\mathbf{S}}^T \widehat{R}^{-1} \widetilde{\mathbf{S}} \xrightarrow{\mathcal{P}} \mathbf{M}_l$ and $n_{l_j}^{1/2} (\widehat{\mathbf{h}} - \mathbf{h}) \xrightarrow{\mathcal{L}} N(\mathbf{0}, \sigma^2 \mathbf{M}_l^{-1})$.

It follows that along $\{n_{l_j}\}_{j=1}^\infty$, $\widehat{\mathbf{h}} \xrightarrow{\mathcal{P}} \mathbf{h}$ as $j \to \infty$. Thus, for any subsequence of $\widehat{\mathbf{h}}$, there exists a further subsequence along which $\widehat{\mathbf{h}} \xrightarrow{\mathcal{P}} \mathbf{h}$. This gives $\widehat{\mathbf{h}} \xrightarrow{\mathcal{P}} \mathbf{h}$ as $n \to \infty$.

We now show the second part. Applying Slutsky's theorem gives that as $j \to \infty$, $\{A(\widetilde{\mathbf{S}}^T \widehat{R}^{-1} \widetilde{\mathbf{S}})^{-1} A^T\}^{-1/2} A(\widehat{\mathbf{h}} - \mathbf{h})$ has an asymptotic Gaussian distribution with mean vector zero and variance–covariance matrix $\sigma^2 \mathbf{I}_k$, which implies that, for $\{n_{l_j}\}_{j=1}^\infty$,

$$(\widehat{\mathbf{h}} - \mathbf{h})^T A^T \{A(\widetilde{\mathbf{S}}^T \widehat{R}^{-1} \widetilde{\mathbf{S}})^{-1} A^T\}^{-1} A(\widehat{\mathbf{h}} - \mathbf{h}) \xrightarrow{\mathcal{L}} \sigma^2 \chi_k^2 \qquad \text{as } j \to \infty.$$

Applying Lemma A.4, we deduce that $(\widehat{\mathbf{h}}-\mathbf{h})^T A^T \{A(\widetilde{\mathbf{S}}^T \widehat{R}^{-1}\widetilde{\mathbf{S}})^{-1} A^T\}^{-1} A(\widehat{\mathbf{h}}-\mathbf{h}) \xrightarrow{\mathcal{L}} \sigma^2 \chi_k^2$ as $n \to \infty$. □

LEMMA A.8. *Assume Condition* A. *Then* $\widehat{\mathbf{r}}^T \widehat{R}^{-1} \widehat{\mathbf{r}}/n \xrightarrow{\mathcal{P}} \sigma^2$.

PROOF. By the definition of $\widehat{\mathbf{r}}$,

$$\widehat{\mathbf{r}} = \widetilde{\mathbf{y}} - \widetilde{\mathbf{S}}\widehat{\mathbf{h}} = (\mathbf{I} - S_d)(\mathbf{y} - \mathbf{S}\widehat{\mathbf{h}}) = \widetilde{\mathbf{S}}(\mathbf{h} - \widehat{\mathbf{h}}) + \widetilde{\mathbf{d}} + \widetilde{\boldsymbol{\epsilon}}. \tag{A.19}$$

Notice that $n^{-1}\widehat{\mathbf{r}}^T \widehat{R}^{-1} \widehat{\mathbf{r}} = n^{-1}\widehat{\mathbf{r}}^T R^{-1} \widehat{\mathbf{r}} + n^{-1}\widehat{\mathbf{r}}^T (\widehat{R}^{-1} - R^{-1})\widehat{\mathbf{r}}$, in which

$$n^{-1}\widehat{\mathbf{r}}^T R^{-1} \widehat{\mathbf{r}} = n^{-1}\|R^{-1/2}\{\widetilde{\mathbf{S}}(\mathbf{h} - \widehat{\mathbf{h}}) + \widetilde{\mathbf{d}} + \widetilde{\boldsymbol{\epsilon}}\}\|^2$$
$$= n^{-1}\|R^{-1/2}\{\widetilde{\mathbf{S}}(\mathbf{h} - \widehat{\mathbf{h}}) + \widetilde{\mathbf{d}} - S_d \boldsymbol{\epsilon}\} + R^{-1/2}\boldsymbol{\epsilon}\|^2$$
$$\equiv I_1 + 2I_2 + I_3,$$

where $I_1 = n^{-1}\|R^{-1/2}\{\widetilde{\mathbf{S}}(\mathbf{h} - \widehat{\mathbf{h}}) + \widetilde{\mathbf{d}} - S_d\boldsymbol{\epsilon}\}\|^2$, $I_2 = n^{-1}\{\widetilde{\mathbf{S}}(\mathbf{h} - \widehat{\mathbf{h}}) + \widetilde{\mathbf{d}} - S_d\boldsymbol{\epsilon}\}^T R^{-1}\boldsymbol{\epsilon}$ and $I_3 = n^{-1}\|R^{-1/2}\boldsymbol{\epsilon}\|^2$. The proof will be completed by showing that $I_1 = o_P(1)$, $I_2 = o_P(1)$, $I_3 = \sigma^2 + o_P(1)$ and $n^{-1}\widehat{\mathbf{r}}^T (\widehat{R}^{-1} - R^{-1})\widehat{\mathbf{r}} = o_P(1)$.

First, consider $I_1$. Note that

$$I_1 = n^{-1}\|R^{-1/2}\{\widetilde{\mathbf{S}}(\mathbf{h} - \widehat{\mathbf{h}}) + \widetilde{\mathbf{d}} - S_d\boldsymbol{\epsilon}\}\|^2$$
$$\leq 3n^{-1}\|R^{-1/2}\widetilde{\mathbf{S}}(\mathbf{h} - \widehat{\mathbf{h}})\|^2 + 3n^{-1}\|R^{-1/2}\widetilde{\mathbf{d}}\|^2 + 3n^{-1}\|R^{-1/2}S_d\boldsymbol{\epsilon}\|^2.$$

The first term is $o_P(1)$ by Lemma A.6 and Corollary A.2, the second and third terms are both $o_P(1)$ by (A.8), (A.15) and similar derivations for (A.7). Thus, $I_1 = o_P(1)$.

Second, consider $I_3 = n^{-1}\sum_{i=1}^n \sum_{j=1}^n \epsilon(t_i)\epsilon(t_j)V(i,j)$. Then,

$$E(I_3) = n^{-1}E(\boldsymbol{\epsilon}^T R^{-1}\boldsymbol{\epsilon}) = n^{-1}\text{trace}\{E(\boldsymbol{\epsilon}\boldsymbol{\epsilon}^T)R^{-1}\} = \sigma^2,$$

$$E(I_3^2) = n^{-2}\sum_{k_1=1}^n \sum_{k_2=1}^n \sum_{k_3=1}^n \sum_{k_4=1}^n e(k_1,k_2,k_3,k_4)V(k_1,k_2)V(k_3,k_4),$$

where $e(k_1,k_2,k_3,k_4) = E\{\epsilon(t_{k_1})\epsilon(t_{k_2})\epsilon(t_{k_3})\epsilon(t_{k_4})\}$. There are only four possible cases in which $e(k_1,k_2,k_3,k_4)$ is nonzero. In Case 1, for any $i \in \{k_1, k_2, k_3, k_4\}$, there exists $j \in \{k_1, k_2, k_3, k_4\}$, $i \neq j$, such that $|i - j| \leq g$. Then,

$$\left| n^{-2}\sum\sum\sum\sum_{\text{Case 1}} e(k_1,k_2,k_3,k_4)V(k_1,k_2)V(k_3,k_4) \right|$$
$$\leq n^{-2}\sum_{k_1=1}^n \sum_{k_2:|k_2-k_1|\leq g} \sum_{k_3:|k_3-k_1|\leq g} \sum_{k_4:|k_4-k_1|\leq g} E[\{\epsilon(t_1)\}^4]$$
$$\times |V(k_1,k_2)| \cdot |V(k_3,k_4)|$$
$$\leq n^{-2}C^2.$$

In Case 2, $|k_1 - k_2| \leq g$ and $|k_3 - k_4| \leq g$, but for any $i \in \{k_1, k_2\}$ and $j \in \{k_3, k_4\}$, $|i - j| > g$. Then,

$$n^{-2} \sum \sum_{\text{Case 2}} \sum \sum e(k_1, k_2, k_3, k_4) V(k_1, k_2) V(k_3, k_4)$$

$$= n^{-2} \sigma^4 \sum_{k_1, k_2} \sum \rho(|k_2 - k_1|) V(k_1, k_2) \sum_{k_3, k_4} \sum \rho(|k_3 - k_4|) V(k_3, k_4)$$

$$- n^{-2} \sigma^4 \sum \sum_{\text{Case 1}} \sum \sum \rho(|k_2 - k_1|) \rho(|k_3 - k_4|) V(k_1, k_2) V(k_3, k_4)$$

$$= \{E(I_3)\}^2 - O(n^{-2}).$$

In Case 3, $|k_1 - k_3| \leq g$ and $|k_2 - k_4| \leq g$, but for any $i \in \{k_1, k_3\}$ and $j \in \{k_2, k_4\}$, $|i - j| > g$. Then,

$$\left| n^{-2} \sum \sum_{\text{Case 3}} \sum \sum e(k_1, k_2, k_3, k_4) V(k_1, k_2) V(k_3, k_4) \right|$$

$$\leq n^{-2} \sigma^4 \sum_{k_1=1}^{n} \sum_{k_3 : |k_3 - k_1| \leq g} \sum_{k_2=1}^{n} \sum_{k_4 : |k_4 - k_2| \leq g} \rho(|k_3 - k_1|) \rho(|k_4 - k_2|)$$

$$\times |V(k_1, k_2)| \cdot |V(k_3, k_4)|$$

$$\leq n^{-1} C^2.$$

A similar result holds for Case 4, where $|k_1 - k_4| \leq g$ and $|k_2 - k_3| \leq g$, but for any $i \in \{k_1, k_4\}$ and $j \in \{k_2, k_3\}$, $|i - j| > g$. Combining Cases 1–4, $E(I_3^2) \to \{E(I_3)\}^2$, which leads to $I_3 \xrightarrow{\mathcal{P}} E(I_3) = \sigma^2$.

Third, consider $I_2$. By the Cauchy–Schwarz inequality, $I_2 = o_P(1)$.

Fourth, to deduce $n^{-1} \hat{\mathbf{r}}^T (\widehat{R}^{-1} - R^{-1}) \hat{\mathbf{r}} = o_P(1)$, it is sufficient to show that

(A.20) $$n^{-1} \{\widetilde{\mathbf{S}}(\mathbf{h} - \widehat{\mathbf{h}})\}^T (\widehat{R}^{-1} - R^{-1}) \widetilde{\mathbf{S}}(\mathbf{h} - \widehat{\mathbf{h}}) = o_P(1),$$

(A.21) $$n^{-1} \widetilde{\mathbf{d}}^T (\widehat{R}^{-1} - R^{-1}) \widetilde{\mathbf{d}} = o_P(1),$$

(A.22) $$n^{-1} \widetilde{\boldsymbol{\epsilon}}^T (\widehat{R}^{-1} - R^{-1}) \widetilde{\boldsymbol{\epsilon}} = o_P(1).$$

It is easy to see that (A.20) follows from (A.16) and Corollary A.2, whereas (A.21)–(A.22) are obtained by similar arguments for (A.16)–(A.18). □

LEMMA A.9. *Assume Condition A. Then:*

1. $n^{-1/2} \widetilde{\mathbf{S}}^T \widehat{R}^{-1} \widetilde{\widehat{\mathbf{d}}} = o_P(1)$;
2. $(A \widehat{\mathbf{h}}_{\text{bc}} - A\mathbf{h})^T \{A (\widetilde{\mathbf{S}}^T \widehat{R}^{-1} \widetilde{\mathbf{S}})^{-1} A^T\}^{-1} (A \widehat{\mathbf{h}}_{\text{bc}} - A\mathbf{h}) \xrightarrow{\mathcal{L}} \sigma^2 \chi_k^2$;
3. $n^{-1} \widehat{\mathbf{r}}_{\text{bc}}^T \widehat{R}^{-1} \widehat{\mathbf{r}}_{\text{bc}} \xrightarrow{\mathcal{P}} \sigma^2$.

PROOF. To show the first part, note that $\widetilde{\widehat{\mathbf{d}}} = (\mathbf{I} - S_d)\widehat{\mathbf{d}} = (\mathbf{I} - S_d)S_d(\mathbf{y} - \mathbf{S}\widehat{\mathbf{h}}) = S_d(\mathbf{I} - S_d)(\mathbf{y} - \mathbf{S}\widehat{\mathbf{h}}) = S_d\widehat{\mathbf{r}}$. Thus from (A.19),

$$\text{(A.23)} \quad \widetilde{\widehat{\mathbf{d}}} = (\mathbf{I} - S_d)S_d\mathbf{S}(\mathbf{h} - \widehat{\mathbf{h}}) + S_d\widetilde{\mathbf{d}} + (\mathbf{I} - S_d)S_d\boldsymbol{\epsilon},$$

in which

$$(\mathbf{I} - S_d)S_d\boldsymbol{\xi}_{j,\ell} = (\mathbf{I} - S_d)S_d[\{\boldsymbol{\xi}_{j,\ell} - E(\boldsymbol{\xi}_{j,\ell})\} - \{p_j\mathbf{1} - E(\boldsymbol{\xi}_{j,\ell})\}]$$
$$= (\mathbf{I} - S_d)[S_d\{\boldsymbol{\xi}_{j,\ell} - E(\boldsymbol{\xi}_{j,\ell})\} - S_d(p_j\mathbf{1}_{j-1}^T, \mathbf{0}_{n-j+1}^T)^T].$$

Note that $n^{-1/2}\widetilde{\mathbf{S}}^T R^{-1}\widetilde{\widehat{\mathbf{d}}} = I_1 + I_2 + I_3$, where $I_1 = n^{-1/2}\widetilde{\mathbf{S}}^T R^{-1}(\mathbf{I} - S_d)S_d\mathbf{S}(\mathbf{h} - \widehat{\mathbf{h}})$, $I_2 = n^{-1/2}\widetilde{\mathbf{S}}^T R^{-1}S_d\widetilde{\mathbf{d}}$ and $I_3 = n^{-1/2}\widetilde{\mathbf{S}}^T R^{-1}(\mathbf{I} - S_d)S_d\boldsymbol{\epsilon}$. We now show that each term is $o_P(1)$. For $I_1$, from Lemma A.7, we have $n^{1/2}(\mathbf{h} - \widehat{\mathbf{h}}) = O_P(1)$, thus we only need to consider the matrix $n^{-1}\widetilde{\mathbf{S}}^T R^{-1}(\mathbf{I} - S_d)S_d\mathbf{S}$. For its block matrix $n^{-1}\widetilde{\mathbf{S}}_{j_1}^T R^{-1}(\mathbf{I} - S_d)S_d\mathbf{S}_{j_2}$, the $(\ell_1, \ell_2)$th entry satisfies

$$|n^{-1}\boldsymbol{e}_{\ell_1}^T \widetilde{\mathbf{S}}_{j_1}^T R^{-1}(\mathbf{I} - S_d)S_d\mathbf{S}_{j_2}\boldsymbol{e}_{\ell_2}|$$
$$\leq n^{-1/2}\|R^{-1/2}\widetilde{\mathbf{S}}_{j_1}\boldsymbol{e}_{\ell_1}\|n^{-1/2}\|R^{-1/2}(\mathbf{I} - S_d)S_d\boldsymbol{\xi}_{j_2,\ell_2}\|$$
$$\leq I_{11}(I_{12} + I_{13}),$$

where $I_{11} = n^{-1/2}\|R^{-1/2}\widetilde{\mathbf{S}}_{j_1}\boldsymbol{e}_{\ell_1}\|$, $I_{12} = n^{-1/2}\|R^{-1/2}(\mathbf{I} - S_d)S_d\{\boldsymbol{\xi}_{j_2,\ell_2} - E(\boldsymbol{\xi}_{j_2,\ell_2})\}\|$ and $I_{13} = n^{-1/2}\|R^{-1/2}(\mathbf{I} - S_d)S_d(p_{j_2}\mathbf{1}_{j_2-1}^T, \mathbf{0}_{n-j_2+1}^T)^T\|$. Then by Lemma A.6, $I_{11} = O_P(1)$. By (A.8), $I_{12} = o(1)$ and, similarly, $I_{13} = o(1)$. Thus $I_1 = o_P(1)$. For $I_2$, using the same procedures as in Lemma A.7 for proving $J_1 = o_P(1)$, we can show $I_2 = o_P(1)$. For $I_3$, using the same procedures as in Lemma A.7 for proving $J_{22} = o_P(1)$, we obtain $I_3 = o_P(1)$. Thus, $n^{-1/2}\widetilde{\mathbf{S}}^T R^{-1}\widetilde{\widehat{\mathbf{d}}} = o_P(1)$. It remains to show that $n^{-1/2}\widetilde{\mathbf{S}}^T(\widehat{R}^{-1} - R^{-1})\widetilde{\widehat{\mathbf{d}}} = o_P(1)$, whose proof is similar to that of (A.17).

To show the second part, recall that $\widehat{\mathbf{h}}_{\text{bc}} = \widehat{\mathbf{h}} - (n^{-1}\widetilde{\mathbf{S}}^T \widehat{R}^{-1}\widetilde{\mathbf{S}})^{-1}(n^{-1}\widetilde{\mathbf{S}}^T \widehat{R}^{-1}\widetilde{\widehat{\mathbf{d}}})$. Using the first part together with Lemma 6 and Corollary 2 leads to the second part.

To show the third part, note that $\widehat{\mathbf{r}}_{\text{bc}} = \widehat{\mathbf{r}} - \widetilde{\widehat{\mathbf{d}}}$ and $n^{-1}\widehat{\mathbf{r}}_{\text{bc}}^T \widehat{R}^{-1}\widehat{\mathbf{r}}_{\text{bc}} = n^{-1}\widehat{\mathbf{r}}_{\text{bc}}^T \times R^{-1}\widehat{\mathbf{r}}_{\text{bc}} + n^{-1}\widehat{\mathbf{r}}_{\text{bc}}^T(\widehat{R}^{-1} - R^{-1})\widehat{\mathbf{r}}_{\text{bc}}$, in which

$$n^{-1}\widehat{\mathbf{r}}_{\text{bc}}^T R^{-1}\widehat{\mathbf{r}}_{\text{bc}} = n^{-1}\|R^{-1/2}(\widehat{\mathbf{r}} - \widetilde{\widehat{\mathbf{d}}})\|^2 \equiv J_1 - 2J_2 + J_3,$$

where $J_1 = n^{-1}\|R^{-1/2}\widehat{\mathbf{r}}\|^2$, $J_2 = 2n^{-1}\widehat{\mathbf{r}}^T R^{-1}\widetilde{\widehat{\mathbf{d}}}$ and $J_3 = n^{-1}\|R^{-1/2}\widetilde{\widehat{\mathbf{d}}}\|^2$. From Lemma A.8, $J_1 = \sigma^2 + o_P(1)$. From (A.23),

$$J_3 = n^{-1}\|R^{-1/2}\{(\mathbf{I} - S_d)S_d\mathbf{S}(\mathbf{h} - \widehat{\mathbf{h}}) + S_d\widetilde{\mathbf{d}} + (\mathbf{I} - S_d)S_d\boldsymbol{\epsilon}\}\|^2$$
$$\leq 3n^{-1}\{\|R^{-1/2}(\mathbf{I} - S_d)S_d\mathbf{S}(\mathbf{h} - \widehat{\mathbf{h}})\|^2$$
$$+ \|R^{-1/2}S_d\widetilde{\mathbf{d}}\|^2 + \|R^{-1/2}(\mathbf{I} - S_d)S_d\boldsymbol{\epsilon}\|^2\}.$$

32C. ZHANG AND T. YU32C. ZHANG AND T. YUUsing similar proofs for the numerator, we obtain $J_3 = o_P(1)$. By the Cauchy–Schwarz inequality, $J_2 = o_P(1)$. Thus, $n^{-1}\widehat{\mathbf{r}}_{\mathrm{bc}}^T R^{-1} \widehat{\mathbf{r}}_{\mathrm{bc}} \xrightarrow{\mathcal{P}} \sigma^2$. To show $n^{-1}\widehat{\mathbf{r}}_{\mathrm{bc}}^T \times (\widehat{R}^{-1} - R^{-1})\widehat{\mathbf{r}}_{\mathrm{bc}} = o_P(1)$, it is sufficient to show that

$$n^{-1}\widehat{\mathbf{r}}^T(\widehat{R}^{-1} - R^{-1})\widehat{\mathbf{r}} = o_P(1), \tag{A.24}$$

$$n^{-1}\widetilde{\widehat{\mathbf{d}}}^T(\widehat{R}^{-1} - R^{-1})\widetilde{\widehat{\mathbf{d}}} = o_P(1), \tag{A.25}$$

where (A.24) directly follows from the fourth step of the proof for Lemma A.8 and (A.25) uses similar proofs for (A.21)–(A.22). □

**A.1. Proof of Theorem 4.1.** From Corollary A.2, under $H_0$ in (4.1), the numerator of $\mathbb{K}$ converges in distribution to $\sigma^2 \chi_k^2$. This, combined with Lemma A.8, gives the desired result for $\mathbb{K}$.

**A.2. Proof of Theorem 4.2.** Under $H_0$ in (4.1), the second and third parts of Lemma A.9 complete the proof for $\mathbb{K}_{\mathrm{bc}}$.

**A.3. Proof of Theorem 4.3.** The numerator of $\mathbb{K}$ can be decomposed into three additive terms:

$$I_1 = (A\widehat{\mathbf{h}} - A\mathbf{h})^T \{A(\widetilde{\mathbf{S}}^T \widehat{R}^{-1} \widetilde{\mathbf{S}})^{-1} A^T\}^{-1} (A\widehat{\mathbf{h}} - A\mathbf{h});$$

$$I_2 = 2n(A\mathbf{h})^T \{A(n^{-1}\widetilde{\mathbf{S}}^T \widehat{R}^{-1} \widetilde{\mathbf{S}})^{-1} A^T\}^{-1} (A\widehat{\mathbf{h}} - A\mathbf{h});$$

$$I_3 = n(A\mathbf{h})^T \{A(n^{-1}\widetilde{\mathbf{S}}^T \widehat{R}^{-1} \widetilde{\mathbf{S}})^{-1} A^T\}^{-1} (A\mathbf{h}).$$

Notice that $I_1 \xrightarrow{\mathcal{L}} \sigma^2 \chi_k^2$ by the second part of Corollary A.2; $I_3 = n(A\mathbf{h})^T \times (A\mathbf{M}^{-1}A^T)^{-1} A\mathbf{h}\{1 + o_P(1)\}$ by Lemma A.7 and $H_1$ in (4.1); $I_2 = O_P(\sqrt{n})$ by the Cauchy–Schwarz inequality. These, along with Lemma A.8, complete the proof for $\mathbb{K}$. The proof for $\mathbb{K}_{\mathrm{bc}}$ is similar and is hence omitted.

**A.4. Proof of Theorem 4.4.** Following Lemma A.7, under $H_{1n}$ in (4.2), $n^{1/2}A\widehat{\mathbf{h}} \xrightarrow{\mathcal{L}} N(\mathbf{c}, \sigma^2 A\mathbf{M}^{-1}A^T)$. Thus

$$\frac{\{A(n^{-1}\widetilde{\mathbf{S}}^T \widehat{R}^{-1} \widetilde{\mathbf{S}})^{-1} A^T\}^{-1/2} n^{1/2} A\widehat{\mathbf{h}}}{(\widehat{\mathbf{r}}^T \widehat{R}^{-1} \widehat{\mathbf{r}})^{1/2}} \xrightarrow{\mathcal{L}} N((A\mathbf{M}^{-1}A^T)^{-1/2}\mathbf{c}/\sigma, \mathbf{I}_k).$$

This completes the proof for $\mathbb{K}$. Similar proofs for $\mathbb{K}_{\mathrm{bc}}$ are omitted.

**Acknowledgments.** We thank Richard Davidson, Tom Johnstone and Terry Oakes at the Waisman Laboratory for Brain Imaging and Behavior, University of Wisconsin-Madison, for helpful comments and for providing the fMRI data set. The comments of the anonymous referee, the Associate Editor and the Co-Editor, Bernard Silverman, were greatly appreciated.

# REFERENCES

BICKEL, P. J. and LEVINA, E. (2008). Regularized estimation of large covariance matrices. *Ann. Statist.* **36** 199–227.

BICKEL, P. J. and LI, B. (2006). Regularization in statistics (with discussion). *Test* **15** 271–344. MR2273731

BENJAMINI, Y. and HOCHBERG, Y. (1995). Controlling the false discovery rate: A practical and powerful approach to multiple testing. *J. Roy. Statist. Soc. Ser. B* **57** 289–300. MR1325392

BOSQ, D. (1998). *Nonparametric Statistics for Stochastic Processes*: *Estimation and Prediction*, 2nd ed. Springer, Berlin. MR1640691

COX, R. W. (1996). AFNI: Software for analysis and visualization of functional magnetic resonance neuroimages. *Comput. Biomed. Res.* **29** 162–73.

DEMKO, S., MOSS, W. F. and SMITH, P. W. (1984). Decay rates for inverses of band matrices. *Math. Comp.* **43** 491–499. MR0758197

FAHRMEIR, L. and GÖSSL, C. (2002). Semiparametric Bayesian models for human brain mapping. *Statistical Modelling* **2** 235–250. MR1951702

FAN, J. and GIJBELS, I. (1996). *Local Polynomial Modeling and Its Applications*. Chapman and Hall, London. MR1383587

FAN, J. and YAO, Q. (2003). *Nonlinear Time Series*: *Nonparametric and Parametric Methods*. Springer, New York. MR1964455

FAN, J., ZHANG, C. M. and ZHANG, J. (2001). Generalized likelihood ratio statistics and Wilks phenomenon. *Ann. Statist.* **29** 153–193. MR1833962

FRISTON, K. J. ET. AL. (1997). SPM course notes. Available at <http://www.fil.ion.ucl.ac.uk/spm/>.

GENOVESE, C. R. (2000). A Bayesian time-course model for functional magnetic resonance imaging data (with discussion). *J. Amer. Statist. Assoc.* **95** 691–703.

GLOVER, G. H. (1999). Deconvolution of impulse response in event-related BOLD fMRI. *NeuroImage* **9** 416–429.

GOLUB, G. H. and VAN LOAN, C. F. (1996). *Matrix Computations*, 3rd ed. Johns Hopkins Univ. Press, Baltimore, MD. MR1417720

GOUTTE, C., NIELSEN, F. A. and HANSEN, L. K. (2000). Modeling the haemodynamic response in fMRI using smooth FIR filters. *IEEE Trans. Med. Imag.* **19** 1188–1201.

JOSEPHS, O. and HENSON, R. N. A. (1999). Event-related functional magnetic resonance imaging: Modelling, inference and optimization. *Philos. Trans. Roy. Soc.* **354** 1215–1228.

LANGE, N. (1996). Tutorial in biostatistics: Statistical approaches to human brain mapping by functional magnetic resonance imaging. *Statistics in Medicine* **15** 389–428.

LANGE, N. and ZEGER, S. L. (1997). Non-linear Fourier times series analysis for human brain mapping by functional magnetic resonance imaging. *Appl. Statist.* **46** 1–29. MR1452285

LAZAR, N. A, EDDY, W. F., GENOVESE, C. R. and WELLING, J. (2001). Statistical issues in fMRI for brain imaging. *Internat. Statist. Rev.* **69** 105–127.

LU, Y. (2006). Contributions to functional data analysis with biological applications. Ph.D. dissertation, Dept. Statistics, Univ. Wisconsin, Madison.

PURDON, P. L., SOLO, V., WEISSKO, R. M. and BROWN, E. (2001). Locally regularized spatio temporal modeling and model comparison for functional MRI. *NeuroImage* **14** 912–923.

SHAO, J. (2003). *Mathematical Statistics*, 2nd ed. Springer, New York. MR2002723

SILVERMAN, B. W. (1986). *Density Estimation For Statistics and Data Analysis*. Chapman and Hall, London. MR0848134

Smith, S., Jenkinson, M., Woolrich, M., Beckmann, C. F., Behrens, T. E. J., Johansen-Berg, H., Bannister, P. R., De Luca, M., Drobnjak, I. Flitney, D. E., Niazy, R. K., Saunders, J., Vickers, J., Zhang, Y., De Stefano, N., Brady, J. M. and Matthews, P. M. (2004). Advances in functional and structural MR image analysis and implementation as FSL. *NeuroImage* **23** 208–219.

Storey, J. D. (2002). A direct approach to false discovery rates. *J. Roy. Statist. Soc. Ser. B* **64** 479–498. MR1924302

Ward, B. D. (2001). Deconvolution analysis of fMRI time series data. Technical report, Biophysics Research Institute, Medical College of Wisconsin.

Woolrich, M. W., Ripley, B. D., Brady, M. and Smith, S. M. (2001). Temporal autocorrelation in univariate linear modelling of fMRI data. *NeuroImage* **14** 1370–1386.

Worsley, K. J. and Friston, K. J. (1995). Analysis of fMRI time-series revisited-again. *NeuroImage* **2** 173–181.

Worsley, K. J., Liao, C. H., Aston, J., Petre, V., Duncan, G. H., Morales, F. and Evans, A. C. (2002). A general statistical analysis for fMRI data. *NeuroImage* **15** 1–15.

Zhang, C. M. (2003). Calibrating the degrees of freedom for automatic data smoothing and effective curve checking. *J. Amer. Statist. Assoc.* **98** 609–628. MR2011675

Zhang, C. M. and Dette, H. (2004). A power comparison between nonparametric regression tests. *Statist. Probab. Lett.* **66** 289–301. MR2045474

Zhang, C. M., Lu, Y., Johnstone, T., Oakes, T. and Davidson, R. (2006). Efficient modeling and inference for event-related fMRI data. Technical report #1125, Dept. Statistics, Univ. Wisconsin-Madison.

Department of Statistics
University of Wisconsin
1300 University Avenue
Madison, Wisconsin 53706
USA
E-mail: cmzhang@stat.wisc.edu
yutao@stat.wisc.edu